\documentclass[10pt]{article}
\usepackage{latexsym,bm}
\usepackage[centertags]{amsmath}
\usepackage{amsfonts}
\usepackage{graphics}
\usepackage{graphicx}
\usepackage{amsthm}
\usepackage{hyperref}
\usepackage{CJK}
\usepackage{times}
\usepackage{mathrsfs}
\usepackage{color}
\usepackage{txfonts}
\def\qed{\hfill ~\qedsymbol}

\newtheorem{prop}{Proposition}[section]
\newtheorem{thm}[prop]{Theorem}
\newtheorem{defn}{Definition}[section]

\newtheorem{lem}[prop]{Lemma}

\newtheorem*{rem}{Remark}

\numberwithin{equation}{section}
\numberwithin{equation}{subsection}

\begin{document}

\title{A mean value formula of  sub-$p$-Laplace parabolic equations on the Heisenberg group
\thanks{This work is supported by National Natural Science Foundation of
China (No.11071119).}}
\author{Hairong Liu$^{1\,2}$\thanks{ E-mail: hrliu@njfu.edu.cn}
\quad Xiaoping Yang$^1$ \thanks{E-mail: yangxp@mail.njust.edu.cn}\\
\small  ${}^{1.}$\,School of Science, Nanjing University of Science
$\&$ Technology, Nanjing 210094,
P. R. China\\
\small  ${}^{2.}$\,School of Science, Nanjing Forestry University,
Nanjing 210037, P. R. China\\
\hspace*{\fill}}
\date{}
 \maketitle {\bf Abstract.}
 We derive two equivalent definitions of the viscosity solutions to the homogeneous sub-$p$-Laplace
 parabolic equations  on the Heisenberg
 group, and characterize the viscosity solutions in terms of an asymptotic mean
value formula, when $1<p\leq \infty$. Moreover, we construct an
example  to show that these formulae do not hold in  non-asymptotic
sense.

{\bf Key Words}\quad Heisenberg  group,   sub-$p$-Laplace parabolic
equation, viscosity solution

{\bf Mathematic Subject Classification.} 35D40, 35K92, 35R03.
\section{Introduction}
Mean value properties for solutions to elliptic and parabolic
partial differential equations are important tools for the study of
their properties. It is well known that  a basic property of
harmonic functions is the mean value property \cite{HL97}. More
precisely,  $u$ is a harmonic function in a domain $\Omega\subset
\mathbb{R}^n$ (that is $u$ satisfies $\Delta u=0$ in $\Omega$) if
and only if $u$ satisfies the mean value formula
\begin{displaymath}
u(x)=\fint_{B_{\varepsilon}(x)}u(y)dy
\end{displaymath}
whenever $B_{\varepsilon}(x)\subset\Omega$  and $\fint_{E}f$ denotes
the average of $f$ over the set $E$.  In addition, an asymptotic
mean value formula holds for some nonlinear cases as well. Manfredi
et al.  \cite{Manfredi2010}  characterized $p$-harmonic functions by
means of asymptotic mean value properties that hold in the so called
viscosity sense (see Definition \ref{def2} below). More precisely,
they proved that the asymptotic mean value formula
\begin{displaymath}
u(x)=\frac{\alpha}{2}\left(\max_{\overline{B}_{\varepsilon}(x)}u
+\min_{\overline{B}_{\varepsilon}(x)}u\right)+\beta\fint_{B_{\varepsilon}(x)}u(y)dy
+o(\varepsilon^2) \hspace{3mm} \mbox{as}\hspace{3mm}
\varepsilon\rightarrow0
\end{displaymath}
 holds  in  the viscosity sense for all
$x\in\Omega$\, if and only if $u$ is a viscosity solution of
\begin{displaymath}
-\Delta_pu=-div\left(|\nabla u|^{p-2}\nabla u\right)=0,
\end{displaymath}
where the constants $\alpha$ and $\beta$ are given by
\begin{displaymath}
\alpha=\frac{p-2}{p+n} \hspace{3mm} \mbox{and} \hspace{3mm}
\beta=\frac{2+n}{p+n}.
\end{displaymath}
The mean value properties of  $p$-Laplace parabolic  equation were
 proved by Manfredi et al.  \cite{manfredi2}. In fact, they
proved that the asymptotic mean value formula
\begin{eqnarray*}
u(t,x)=\frac{\alpha}{2}\fint_{t-\varepsilon^2}^{t}\left(\max_{y\in
\overline{B}_{\varepsilon}(x)}u(s,y)+\min_{y\in
\overline{B}_{\varepsilon}(x)}u(s,y)\right)ds
+\beta\fint_{t-\varepsilon^2}^{t}\fint_{B_{\varepsilon}(x)}u(s,y)dyds+o(\varepsilon^2),
\mbox{as}\hspace{1mm}
\varepsilon\rightarrow0
\end{eqnarray*}
holds for every $(t,x)\in \Omega_{T}=(0,T)\times \Omega$ in the
viscosity sense if and only if $u$ is a viscosity solution of
\begin{equation*}
(n+p)u_t(t,x)=|\nabla u|^{2-p}\Delta_{p}u(t,x).
\end{equation*}
The constants $\alpha$ and $\beta$ are the same as before.

The purpose of  this paper is to extend this result to parabolic
equations on the Heisenberg group  $\mathbb{H}^n$. We recall that
$\mathbb{H}^n$ is the Lie group $(\mathbb{R}^{2n+1},\circ)$ equipped
with the group action
\begin{equation}\label{action}
x^0\circ x=\left(x_1+x_1^0,\cdots,
x_{2n}+x_{2n}^0,x_{2n+1}+x_{2n+1}^0+2\sum\limits_{i=1}^n(x_ix_{n+i}^0-x_i^0x_{n+i})\right),
\end{equation}
for
$x=(x_1,\cdots,x_n,x_{n+1},\cdots,x_{2n},x_{2n+1})=(\overline{x},x_{2n+1})\in
\mathbb{R}^{2n+1}$. It is easy to check that (\ref{action}) does
indeed make $\mathbb{R}^{2n}\times \mathbb{R}$ into a group whose
identity is the origin, and where the inverse is given by
$x^{-1}=-x$. Let us denote by $\delta_{\lambda}$ the Heisenberg
group dilation
\begin{equation}\label{dilation}
\delta_{\lambda}(x_1,\cdots,x_{2n},x_{2n+1})=(\lambda
x_1,\cdots,\lambda x_{2n},\lambda^2 x_{2n+1}),\hspace{2mm}
\lambda>0.
\end{equation}
Then $\mathbb{H}^n= (\mathbb{R}^{2n+1},\circ, \delta_{\lambda})$ is
a homogeneous group.  We denote $Q = 2n + 2$ and call it the
homogeneous dimension of $\mathbb{H}^n$.  For more information on
the Heisenberg group, we refer the reader to the monograph
\cite{BLU07}.

A basis of the Lie algebra of $\mathbb{H}^n$ is
given by
\begin{equation}\label{vector}
\left\{
\begin{array}{lll}
 \displaystyle X_i=\frac{\partial}{\partial x_i}+2x_{n+i}\frac{\partial}{\partial
x_{2n+1}}, \hspace{2mm}i=1,\cdots,n, \\[2mm]
 \displaystyle X_{n+i}=\frac{\partial}{\partial
x_{n+i}}-2x_i\frac{\partial}{\partial
x_{2n+1}}, \hspace{2mm}i=1,\cdots,n, \\[2mm]
 \displaystyle T=\frac{\partial}{\partial x_{2n+1}}.
\end{array}
\right.
\end{equation}
From (\ref{vector}), it is easy to check that  $X_i$ and $X_{n+j}$
satisfy
\begin{displaymath}
[X_i,X_{n+j}]=-4T\delta_{ij},\hspace{2mm}[X_i,X_j]=[X_{n+i},X_{n+j}]=0,\hspace{2mm}
i,j=1,\cdots,n.
\end{displaymath}
Therefore, the vector fields $X_i$, $X_{n+i}$  $(i=1,\cdots,n)$ and
their first order commutators span the whole Lie Algebra.

 Given a
function $u: \mathbb{H}^n\rightarrow R$, we consider the full
gradient of $u$
\begin{equation*}
\nabla u=(X_1u,\cdots, X_{2n}u, Tu)
\end{equation*}
 and the horizontal gradient of $u$
\begin{equation*}
\nabla_{0}u=(X_1u,\cdots,X_{2n}u)
\end{equation*}
 and the symmetrized second
horizontal derivative matrix $(X^2u)^*$
\begin{displaymath}
(X^2u)^*=\frac{1}{2}\left(X_iX_ju+X_jX_iu\right).
\end{displaymath}

 For $x \in \mathbb{H}^n$,
we define the quasi-distance from the origin
\begin{displaymath}
\rho(x)=\left(\left(\sum\limits_{i=1}^n(x_i^2+x_{n+i}^2)\right)^2+x_{2n+1}^2\right)^{\frac{1}{4}}\equiv
\left(|\overline{x}|^4+x_{2n+1}^2\right)^{\frac{1}{4}},
\end{displaymath}
which satisfies $\rho(\delta_{\lambda}(x)) = \lambda \rho(x)$ and
means that $\rho$ is homogeneous of degree one with respect to the
dilation $\delta_{\lambda}$. The associated distance between $x$ and
$x^0$  is defined by
\begin{displaymath}
\rho(x;x^0)=\rho\left((x^0)^{-1}\circ x\right).
\end{displaymath}

In the sequel we let
\begin{displaymath}
B_r=\{x\in \mathbb{H}^n|\rho(x)<r\}, \hspace{2mm}\partial B_r=\{x\in
\mathbb{H}^n|\rho (x)=r\}
\end{displaymath}
and call these sets  a Heisenberg-ball and a sphere centered at the
origin with radius $r$ respectively. Balls and spheres centered at
$x^0$ are defined by left-translation, i.e.
\begin{displaymath}
B_r(x^0)=\{x\in \mathbb{H}^n|\rho(x;x^0)<r\}, \partial B_r(x^0)=\{x\in \mathbb{H}^n|\rho(x;x^0)=r\}.
\end{displaymath}

Introducing the function
\begin{equation}\label{psi}
\psi(x)=|\nabla_0\rho|^2=\frac{|\overline{x}|^2}{\rho(x)^2},
\end{equation}
we define
\begin{equation}\label{volume}
|B_r|=\fint_{B_r}\psi dx \hspace{2mm}\mbox{and}\hspace{2mm}
|\partial B_r|=\frac{d}{dr}|B_r|.
\end{equation}

Gaveau \cite{Gaveau77} proved the following mean value formula for
the sub-Laplace equation
\begin{equation}\label{sub-laplace}
\Delta_{H}u=\sum_{i=1}^{2n}X_iX_iu=0
\end{equation}
on $\mathbb{H}^n$: let $u$ solve the equation $\Delta_{H}u(x)=0$,
then
\begin{equation}\label{p=2}
u(x^0)=\fint_{B_r(x^0)} \psi\left((x^0)^{-1}\circ x\right)u(x) dx,
\end{equation}
where
\begin{displaymath}
\psi\left((x^0)^{-1}\circ
x\right)=\frac{\left|\overline{x}-(\overline{x^0})\right|^2}{\rho\left((x^0)^{-1}\circ(x)\right)^2}.
\end{displaymath}
Recently in \cite{liu}, we characterized sub-$p$-harmonic functions
on $\mathbb{H}^n$  by  asymptotic mean value formulae  in the
viscosity sense. More precisely, we proved that the asymptotic mean
\begin{displaymath}
u(x^0)=\frac{\alpha}{2}\left(\max_{\overline{B}_{\varepsilon}(x^0)}u+\min_{\overline{B}_{\varepsilon}(x^0)}u\right)
+\beta \fint_{B_{\varepsilon}(x^0)}\psi\left((x^0)^{-1}\circ
x\right)u(x)dx+o(\varepsilon^2)
\end{displaymath}
holds as $\varepsilon\rightarrow 0$ for all $x^0\in \Omega$ in the
viscosity sense if and only if $u$ is a viscosity solution of
\begin{equation}\label{sub-p-laplace}
-\Delta_{H}^pu(x)=\sum_{i=1}^{2n}X_i\left(|\nabla_0
u|^{p-2}X_iu\right)(x)=0.
\end{equation}

\vspace{1mm}

In this paper, we study the  parabolic version of the
sub-$p$-Laplace equation on $\mathbb{H}^n$:
\begin{equation*}
u_{t}(t,x)=|\nabla_0u|^{2-p}\Delta_{H}^{p}u(t,x).
\end{equation*}
Recall that for $1<p<\infty$, we have
\begin{equation}\label{operator1}
u_{t}(t,x)=|\nabla_0u|^{2-p}\Delta_{H}^{p}u=(p-2)\Delta_{H}^{\infty}u+\Delta_{H}u,
\end{equation}
where
\begin{equation}\label{operator2}
\Delta_{H}^{\infty}u=|\nabla_0u|^{-2}\left<\left(X^2u\right)^*\nabla_0u,\nabla_0u\right>
=\left|\nabla_0u\right|^{-2}\sum_{i,j=1}^{2n}X_iX_ju\cdot X_iu\cdot
X_ju
\end{equation}
denotes the 1-homogeneous version of sub-infinity Laplace equation
on $\mathbb{H}^n$.

Before proceeding, we would like to mention some motivations related
to  our research. Since H\"{o}rmander¡¯s work \cite{Hormander1967}
the study of partial differential equations of sub-elliptic type
like (\ref{sub-laplace}), \eqref{sub-p-laplace} and
\eqref{operator2} has received a strong impulse, see,
e.g.,\cite{Bieske2004}, \cite{Bieske2006},
\cite{Bony1969},\cite{DM2005}, \cite{D2008},
\cite{domokos},\cite{Folland and Stein}, \cite{Stein1976},
\cite{wang2007} etc. These equations arise in many different
settings: geometric theory of several complex variables, curvature
problems for CR-manifolds, sub-Riemannian geometry, diffusion
processes, control theory, human vision; see, e.g.,\cite{Citti},
\cite{Jerison}. The parabolic counterpart of the operator is also of
great relevance; see, e.g.,
 \cite{Alexopoulos}, \cite{BLU2003}, \cite{huisken}, \cite{montanari}.\\

Let $T>0$ and $\Omega\subset \mathbb{H}^n$ be an open set, and let
$\Omega_{T}= (0,T)\times\Omega$ be a space-time cylinder. Our main
results are the following  theorems corresponding to  $p=2$,
$p=\infty$ and $1<p<\infty$, respectively.
\begin{thm}\label{thm-p=2}
Let $u$ be a smooth function in $\Omega_{T}$. The asymptotic mean
value formula
\begin{equation}\label{vis-p=2}
u(t,x)=\fint_{t-\varepsilon^2}^{t}
\fint_{B_{\varepsilon}(x)}\psi(x^{-1}\circ
y)u(s,y)dyds+o(\varepsilon^2)\hspace{2mm} \mbox{as}\hspace{2mm}
\varepsilon\rightarrow 0
\end{equation}
holds for all $(t,x)\in \Omega_{T}$  if and only if
\begin{equation}\label{equationp=2}
u_t(t,x)=M(n)\Delta_{H}u(t,x)
\end{equation}
in $\Omega_{T},$ where
\begin{equation}\label{constant}
 M(n)=\left\{\begin{array}{lll}\displaystyle
M_o(n)=\frac{2n+2}{2n+4}\cdot\frac{(n!!)^2}{(n+1)!!(n-1)!!}\cdot\frac{\pi}{2}\cdot\frac{1}{2n},\hspace{2mm}\mbox{if}
\hspace{1mm}n \hspace{1mm}\mbox{is odd},\\[4mm]
\displaystyle
M_e(n)=\frac{2n+2}{2n+4}\cdot\frac{(n!!)^2}{(n+1)!!(n-1)!!}\cdot\frac{2}{\pi}\cdot\frac{1}{2n},\hspace{2mm}\mbox{if}
\hspace{1mm}n\hspace{1mm} \mbox{is even}.
\end{array}\right.
\end{equation}
\end{thm}

Next, we study the homogeneous  sub-infinity Laplace parabolic
equation
\begin{equation}\label{parabolic-infinity}
u_t=\Delta_{H}^{\infty}u=|\nabla_0u|^{-2}\left<(X^2u)^*\nabla_0u,\nabla_0u\right>.
\end{equation}
Since the right-hand side of equation (\ref{parabolic-infinity})
cannot be in a divergence form, we are not able to define a
distributional weak solution. However,  there is a standard way to
define viscosity solutions for singular parabolic equations. We
recall this definition  follow Evans and Spruck \cite{ES}, Chen,
Giga and Goto \cite{CGG}, Ohnuma and Sato \cite{OS}, etc. In
addition, the homogenous sub-infinity Laplace equation
\begin{equation*}
|\nabla_0u|^{-2}\left<(X^2u)^*\nabla_0u,\nabla_0u\right>=0
\end{equation*}
 is different from the  inhomogeneous sub-infinity Laplace equation
\begin{equation*}
\left<(X^2u)^*\nabla_0u,\nabla_0u\right>=0,
\end{equation*}
 which was studied by Bieske \cite{bieske2} and Wang \cite{wang2007}. The primary
difficulty arising from the homogenous sub-infinity Laplace equation
will be to modify the theory to cover the possibility that the
spatial horizontal gradient $\nabla_0u$ may vanish.

For a symmetric matrix $A$, we denote its largest and smallest
eigenvalue by $\lambda_{\max}(A)$ and $\lambda_{\min}(A)$,
respectively. We give a definition of viscosity solutions to
equation (\ref{parabolic-infinity}) as follows:
\begin{defn}\label{solution-of-infinity}
A lower semi-continuous function $u:\Omega_{T}\rightarrow
\mathbb{R}\cup\{+\infty\}$ is a viscosity super-solution to
(\ref{parabolic-infinity}) if for every $(t^0,x^0)\in\Omega_{T}$ and
$\phi\in
C_H^2(\Omega_{T})$ satisfy the following \\[1mm]
(i) $u$ is not identically infinity in each component of $\Omega_{T}$,\\[1mm]
(ii) $u(t^0,x^0)=\phi(t^0,x^0)$, and  $u(t,x)>\phi(t,x) \hspace{2mm}\mbox{for} \hspace{2mm}(t,x)\neq (t^0,x^0),\hspace{2mm} $\\[1mm]
then we have at the point $(t^0,x^0)$
\begin{equation}
\left\{
\begin{array}{lll}
\phi_{t}\geq\Delta_{H}^{\infty}\phi\hspace{15mm}
if \hspace{1mm} \nabla_{0}\phi(t^0,x^0)\neq0, \nonumber \\[2mm]
\phi_{t}\geq\lambda_{min}((X^2\phi)^{*})\hspace{3mm}
if \hspace{1mm} \nabla_{0}\phi(t^0,x^0)=0. \nonumber\\
\end{array}
\right.
\end{equation}

A function $u$ is a viscosity sub-solution to
(\ref{parabolic-infinity}) if $-u$ is a viscosity super-solution. A
function $u$ is a viscosity solution if it is both a viscosity
super-solution and a viscosity sub-solution.
\end{defn}

Similarly to the case in \cite{manfredi2}, the asymptotic mean value
formulae hold in a viscosity sense. We recall the following
definition \cite{manfredi2}.
\begin{defn}
\label{def2}
 A continuous function $u$ satisfies
 \begin{equation}\label{vis}
u(t,x)=\frac{\alpha}{2}\fint_{t-\varepsilon^2}^{t}\left(\max_{y\in\overline{B}_{\varepsilon}(x)}u(s,y)
+\min_{y\in\overline{B}_{\varepsilon}(x)}u(s,y)\right)ds
+\beta\fint_{t-\varepsilon^2}^{t}
\fint_{B_{\varepsilon}(x)}\psi(x^{-1}\circ
y)u(s,y)dsdy+o(\varepsilon^2)
\end{equation}
as $\varepsilon\rightarrow 0$ in the viscosity sense if for every
$\phi\in C_{H}^2$ such that $u-\phi$ has a strict minimum at the
point $(x,t)\in \Omega_{T}$ with $u(x,t)=\phi(x,t)$, we have
\begin{align}
\phi(t,x)\geq\frac{\alpha}{2}\fint_{t-\varepsilon^2}^{t}\left(\max_{y\in\overline{B}_{\varepsilon}(x)}\phi(s,y)
+\min_{y\in\overline{B}_{\varepsilon}(x)}\phi(s,y)\right)ds
+\beta\fint_{t-\varepsilon^2}^{t}
\fint_{B_{\varepsilon}(x)}\psi(x^{-1}\circ y)\phi(s,y)
dsdy+o(\varepsilon^2)
\end{align}
as  $\varepsilon\rightarrow 0,$ and analogously when testing from
above.
\end{defn}

Observe that a $C^{2}_{H}$ function ( see Definition \ref{def-con})
$u$ satisfies an equality in the classical sense if and only if it
satisfies in the viscosity sense.

\begin{thm}\label{thm-infinity}
Let $u$ be a continuous function in $\Omega_{T}$. The asymptotic
mean value formula\begin{equation}\label{vis-infinity}
u(t,x)=\frac{1}{2}\fint_{t-\varepsilon^2}^{t}
\left(\max_{y\in\overline{B}_{\varepsilon}(x)}u(s,y)
+\min_{y\in\overline{B}_{\varepsilon}(x)}u(s,y)\right)ds+o(\varepsilon^2)\hspace{2mm}
\mbox{as}\hspace{2mm} \varepsilon\rightarrow 0
\end{equation}
holds for all $(t,x)\in \Omega_{T}$ in the viscosity sense if and
only if $u$ is  a viscosity solution to (\ref{parabolic-infinity}).
\end{thm}

Finally, we combine the above results to obtain an asymptotic mean
value formula of sub-$p$-Laplace parabolic equations. Recalling the
following definition of viscosity solutions.

\begin{defn}\label{def-p}
A lower semi-continuous function $u:\Omega_{T}\rightarrow
R\cup\{+\infty\}$ is a viscosity super-solution to (\ref{operator1})
if for every $(t^0,x^0)\in\Omega_{T}$ and $\phi\in
C_H^2(\Omega_{T})$ satisfy the following \\[1mm]
(i) $u$ is not identically infinity in each component of
$\Omega_{T}$,\\[1mm]
(ii) $u(t^0,x^0)=\phi(t^0,x^0)$, and  $u(t,x)>\phi(t,x)
\hspace{2mm}\mbox{for} \hspace{2mm}(t,x)\neq (t^0,x^0),\hspace{2mm}
$\\[1mm]
then we have at the point $(t^0,x^0)$
\begin{equation}
\left\{
\begin{array}{lll}
\displaystyle
\phi_{t}\geq(p-2)\Delta_{H}^{\infty}\phi+\Delta_{H}\phi\hspace{15mm}
if \hspace{1mm} \nabla_{0}\phi(t^0,x^0)\neq0, \nonumber\\[2mm]
\displaystyle
\phi_{t}\geq\lambda_{min}((p-2)(X^2\phi)^{*})+\Delta_{H}\phi\hspace{3mm}
if \hspace{1mm} \nabla_{0}\phi(t^0,x^0)=0. \nonumber\\
\end{array}
\right.
\end{equation}
A function $u$ is a viscosity sub-solution to (\ref{operator1}) if
$-u$ is a viscosity super-solution. A function $u$ is a viscosity
solution if it is both a viscosity super-solution and a viscosity
sub-solution.
\end{defn}

We derive  an equivalent definition of the above definition of
viscosity solutions by reducing the number of test functions. We
will prove that, in the case $\nabla_{0}\phi(t,x)=0$, we may assume
that $(X^2\phi)^{*}(t,x)=0$, and thus
$\lambda_{min}=\lambda_{max}=0$. Nothing is required if
$\nabla_{0}\phi(t,x)=0$ and $(X^2\phi)^{*}(t,x)\neq0$. Indeed, we
have
\begin{thm}\label{thm2-1}
Suppose $u:\Omega_{T}\rightarrow \mathbb{R}$ is a lower
semi-continuous function with the property that for every
$(t^0,x^0)\in\Omega_{T}$ and $\phi\in C_{H}^2(\Omega_{T})$
satisfying
 \begin{equation*}
 u(t^0,x^0)=\phi(t^0,x^0)\hspace{2mm}\mbox{and} \hspace{2mm} u(t,x)>\phi(t,x) \hspace{2mm}\mbox{for}
 \hspace{2mm}(t,x)\neq (t^0,x^0),
 \end{equation*}
  the following holds:
\begin{equation}
\left\{
\begin{array}{lll}
\displaystyle
\phi_{t}(t^0,x^0)\geq(p-2)\Delta_{H}^{\infty}\phi(t^0,x^0)+\Delta_{H}\phi(t^0,x^0)\hspace{3mm}
if \hspace{1mm} \nabla_0\phi(t^0,x^0)\neq0, \nonumber \\[2mm]
\displaystyle \phi_{t}(t^0,x^0)\geq 0\hspace{35mm} if \hspace{1mm}
\nabla_{0}\phi(t^0,x^0)=0, \hspace{1mm}\mbox{and} \hspace{1mm}
 (X^2\phi)^{*}(t^0,x^0)=0.\nonumber
\end{array}
\right.
\end{equation}
Then $u$ is a viscosity super-solution of (\ref{operator1}). And the
same result holds for the viscosity sub-solution.
\end{thm}

Theorem \ref{thm-p=2}, together with Theorem \ref{thm-infinity},
immediately yields the following asymptotic mean value formula of
sub-$p$-Laplace parabolic equations.
\begin{thm}\label{thm1}
Let $1<p\leq\infty$ and  $u$ be a continuous function in
$\Omega_{T}$. The asymptotic expansion
\begin{align*}
u(t,x)=\frac{\alpha}{2}\fint_{t-\varepsilon^2}^{t}\left(\max_{y\in\overline{B}_{\varepsilon}(x)}u(s,y)
+\min_{y\in\overline{B}_{\varepsilon}(x)}u(s,y)\right)ds
+\beta\fint_{t-\varepsilon^2}^{t}
\fint_{B_{\varepsilon}(x)}\psi(x^{-1}\circ
y)u(s,y)dsdy+o(\varepsilon^2)\hspace{1mm} \mbox{as}\hspace{1mm}
\varepsilon\rightarrow 0
\end{align*}
holds for all $(t,x)\in \Omega_{T}$ in the viscosity sense if and
only if $u$ is  a viscosity solution to
\begin{displaymath}
u_t(t,x)=\frac{M(n)}{M(n)(p-2)+1}\left|\nabla_0u\right|^{2-p}\Delta_{H}^{p}u(t,x),
\end{displaymath}
where $M(n)$ is as in \eqref{constant}, and $\alpha$ and $\beta$
satisfy
\begin{equation}\label{constant1}
\left\{
\begin{array}{lll}
\beta M(n)(p-2)=\alpha, \\
 \alpha+\beta=1.
\end{array}
\right.
\end{equation}
\end{thm}

\begin{rem}
If $p=2$, then $\alpha=0$ and $\beta=1$, and if $p=\infty$, then
$\alpha=1$ and $\beta=0$.
\end{rem}

The rest of the paper is organized as follows. In Section 2 we
collect some definitions and results about sub-parabolic jets on
$\mathbb{H}^n$. Using the polar coordinates on $\mathbb{H}^n$, we
compute  some integrals.  By twisting the Euclidean jets to
sub-parabolic jets and using the Crandall-Ishii-Lions maximum
principle, Theorem \ref{thm2-1} is proved in Section 3. In Section 4
we prove asymptotic  mean value formulae of sub-heat equations,
sub-infinity Laplace parabolic equation and sub-$p$-Laplace
parabolic equation, respectively.  An example is constructed to show
that these formulae do not hold in non-asymptotic sense.

\section{ Sub-parabolic  jets and   polar
coordinates on $\mathbb{H}^n$}

In this section, we collect some definitions and results about
sub-parabolic jets on $\mathbb{H}^n$,  and recall the polar
coordinates on $\mathbb{H}^n$.

\begin{defn}(\cite{FE82})\label{def-con}
Let $f: \mathbb{H}^n\rightarrow \mathbb{R},$ we say that $f\in
C_{H}^1$,
 if $X_if$ exists and is continuous at each point of $\mathbb{H}^n$, for
every $i=1,\cdots,2n$.  Moreover, for any nonnegative integer $m$,
we say that $f\in C_{H}^m$, if $X^{\alpha}f$ exists and is
continuous at each point of $\mathbb{H}^n$, for every horizontal
derivation $X^{\alpha}=X_{i_1}^{\alpha_1}\cdots
X_{i_n}^{\alpha_n}X_{i_{n+1}}^{\alpha_{n+1}}\cdots
X_{i_{2n}}^{\alpha_{2n}}$ with
$0\leq|\alpha|=\alpha_1+\cdots+\alpha_{2n}\leq m$. \label{def}
\end{defn}

Let $\mathcal{S}^{n}$ be the the set of all real $n\times n$
symmetric matrixes, we introduce definitions about  sub-parabolic
jets on $\mathbb{H}^n$,  which are  natural extensions of
sub-elliptic jets \cite{bieske2}.
\begin{defn}\label{heisenberg-jet}
Let $u:\Omega_{T}\rightarrow \mathbb{R}$ be an upper-semicontinuous
function. The second order sub-parabolic super-jet of $u$ at
$(t^0,x^0)$ is defined as
\begin{align*}
\mathcal{J}^{2,+}\left(u,(t^0,x^0)\right)&=\Big\{(a,p,Y)\in
\mathbb{R}\times \mathbb{R}^{2n+1}\times \mathcal {S}^{2n}
\hspace{2mm}\mbox{such
that}\hspace{2mm}\\[1mm]
&u\left(t,x\right)\leq u\left(t^0,x^0\right)+a\left(t-t^0\right)+\left<p,\left((x^{0})^{-1}\circ x\right)\right>\\[1mm]
&+\frac{1}{2}\left<Y\overline{\left((x^{0})^{-1}\circ
x\right)},\overline{\left((x^{0})^{-1}\circ x\right)}\right>
+c\left(|t-t^0|+\rho^2\left((x^{0})^{-1}\circ x\right)\right)\Big\}.
\end{align*}
Similarly, for a  lower-semicontinuous function $u$, we define the
second order sub-parabolic sub-jet
\begin{align*}
\mathcal{J}^{2,-}\left(u,(t^0,x^0)\right)&=\Big\{(a,p,Y)\in
\mathbb{R}\times \mathbb{R}^{2n+1}\times \mathcal
{S}^{2n}\hspace{2mm}\mbox{such
that}\hspace{2mm}\\[1mm]
&u\left(t,x\right)\geq u\left(t^0,x^0\right)+a\left(t-t^0\right)+\left<p,\left((x^{0})^{-1}\circ x\right)\right>\\[1mm]
&+\frac{1}{2}\left<Y\overline{\left((x^{0})^{-1}\circ
x\right)},\overline{\left((x^{0})^{-1}\circ x\right)}\right>
+c\left(|t-t^0|+\rho^2\left((x^{0})^{-1}\circ x\right)\right)\Big\}.
\end{align*}
The closures of the jets is defined in the obvious way:
\begin{align*}
\overline{\mathcal{J}^{2,+}}\left(u,(t^0,x^0)\right)=\Big\{(a,p,Y)\in
\mathbb{R}\times \mathbb{R}^{2n+1}\times \mathcal {S}^{2n}: \exists
 \left(t^{\alpha},x^{\alpha},a^{\alpha},p^{\alpha},Y^{\alpha}\right)\in
\Omega_{T}\times\mathbb{R}\times \mathbb{R}^{2n+1}\times \mathcal {S}^{2n}\\
\mbox{such that}\hspace{1mm}
(a^{\alpha},p^{\alpha},Y^{\alpha})\in\mathcal{J}^{2,+}\left(u,(t^{\alpha},x^{\alpha})\right)\hspace{1mm}
\hbox{and}\hspace{1mm}
(t^{\alpha},x^{\alpha},a^{\alpha},p^{\alpha},Y^{\alpha})\rightarrow
(t^{0},x^{0},a,p,Y)\Big\},
\end{align*}
and  similarly for $\overline{\mathcal{J}^{2,-}}$.
\end{defn}

 The following proposition characterizes the sub-parabolic jets in
terms of test functions that touch from above or below. This
proposition is an natural extension of the sub-elliptic case
\cite{bieske2}.
\begin{prop}
Define the set
\begin{align*}
K^{2,+}\left(u,(t^0,x^0)\right)=\Big\{\left(\phi_{t}(t^0,x^0),\nabla\phi(t^0,x^0),
(X^2\phi)^{*}(t^0,x^0)\right): \phi\in C^2_{H}(\Omega_{T})\hspace{2mm} \mbox{and}\\
u-\phi\hspace{1mm} \mbox{has a strict  maximum at}
\hspace{1mm}(t^0,x^0)\Big\},
\end{align*}
and
\begin{eqnarray*}
K^{2,-}\left(u,(t^0,x^0)\right)=\Big\{\left(\phi_{t}(t^0,x^0),\nabla\phi(t^0,x^0),
(X^2\phi)^{*}(t^0,x^0)\right): \phi\in C^2_{H}(\Omega_{T})\hspace{2mm} \mbox{and}\\[1mm]
u-\phi\hspace{2mm} \mbox{has a strict  minimum at}
\hspace{2mm}(t^0,x^0)\Big\}.
\end{eqnarray*}
 Then, we have
\begin{equation*}
\mathcal{J}^{2,+}\left(u,(t^0,x^0)\right)=K^{2,+}\left(u,(t^0,x^0)\right),
\end{equation*}
and
\begin{equation*}
\mathcal{J}^{2,-}\left(u,(t^0,x^0)\right)=K^{2,-}\left(u,(t^0,x^0)\right).
\end{equation*}
\end{prop}


At the last of this section, we  recall  polar coordinates on
$\mathbb{H}^n$, which were introduced for $\mathbb{H}^1$ by
\cite{Greiner1980} and then extended by Dunkl \cite{Dunkl} to
$\mathbb{H}^n$.
 Let
\begin{equation}\label{polar2}
\left\{ \begin{array}{lll}
      \displaystyle x_1=\rho\sin^{1/2}\phi \sin\theta_1 \cdots\sin\theta_{2n-2} \sin\theta_{2n-1},   \\[1mm]
      \displaystyle  x_{n+1}= \rho\sin^{1/2}\phi \sin\theta_1 \cdots\sin\theta_{2n-2} \cos\theta_{2n-1},\\[1mm]
   x_2=\rho\sin^{1/2}\phi \sin\theta_1
     \cdots\sin\theta_{2n-3}\cos\theta_{2n-2},\\[1mm]
       \displaystyle x_{n+2}=\rho\sin^{1/2}\phi \sin\theta_1
       \cdots\sin\theta_{2n-4}\cos\theta_{2n-3},\\[1mm]
           \displaystyle\vdots\\
          \displaystyle x_n=\rho\sin^{1/2}\phi \sin\theta_1\cos\theta_2,\\
         \displaystyle x_{2n}=\rho\sin^{1/2}\phi \cos\theta_1,\\[1mm]
          \displaystyle x_{2n+1}=\rho^2\cos\phi.
         \end{array}  \right.
\end{equation}
Here $0\leq\phi<\pi$, $0\leq\theta_{i}<\pi$,  $i=1,\cdots,2n-2$ and
$0\leq\theta_{2n-1}<2\pi$.  Let
$r=|\overline{x}|=\left(\sum_{i=1}^{2n}x_i^2\right)^{1/2}$, from
(\ref{polar2}) we get
\begin{equation}\label{r}
r=|\overline{x}|=\rho \sin^{1/2}\phi.
\end{equation}
By the usual spherical coordinates in $\mathbb{R}^{2n}$, we have
\begin{displaymath}
d\overline{x}=r^{2n-1}dr d\omega,
\end{displaymath}
where $d\omega$ denotes the Lebesgue measure on $S^{2n-1}$. From
(\ref{polar2}) and (\ref{r}) we have
\begin{displaymath}
drdt=\rho^2sin^{-1/2}\phi d\rho d\phi.
\end{displaymath}
Therefore, the Jacobi of (\ref{polar2}) is
\begin{align}\label{jacobi}
dx&=\rho^{2n+1}(\sin\phi)^{n-1} d\rho d\phi dw\nonumber\\
&=\rho^{2n+1}(\sin\phi)^{n-1}\sin^{2n-2}\theta_1\cdots\sin\theta_{2n-2}d\rho
d\phi d\theta_1 \cdots d\theta_{2n-1}.
\end{align}
Using the the polar coordinates on $\mathbb{H}^n$, we calculate to
obtain the following Lemma.
\begin{lem}\label{lem2}
\begin{equation}\label{1}
\int_{B_{\varepsilon}(x)}\psi(x^{-1}\circ y)(y_i-x_i)dy=0,
\hspace{3mm} i=1,\cdots,2n,\nonumber
\end{equation}
\begin{equation}\label{2}
\int_{B_{\varepsilon}(x)}\psi(x^{-1}\circ
y)\left(y_{2n+1}-x_{2n+1}+2\sum_{i=1}^{n}(x_iy_{n+i}-x_{n+i}y_i)\right)dy=0,\nonumber
\end{equation}
and
\begin{equation}\label{3}
\int_{B_{\varepsilon}(x)}\psi(x^{-1}\circ y)(y_i-x_i)\cdot
(y_j-x_j)dy=0\hspace{2mm} \mbox{for} \hspace{2mm}
\hspace{1mm}i,j=1,\cdots,2n, \hspace{1mm} i\neq j.\nonumber
\end{equation}
For every $i=1,\cdots,2n$,  if $n$ is even,
\begin{equation}\label{M_e(n)}
\fint_{B_{\varepsilon}(x)}\psi(x^{-1}\circ y)(y_i-x_i)^2dy
=\frac{2n+2}{2n+4}\cdot\frac{(n!!)^2}{(n+1)!!(n-1)!!}\cdot\frac{1}{2n}\cdot\frac{2}{\pi}\varepsilon^2\equiv
M_e(n)\varepsilon^2,\nonumber
\end{equation}
and if $n$ is odd,
\begin{equation}\label{M_o(n)}
\fint_{B_{\varepsilon}(x)}\psi(x^{-1}\circ y)(y_i-x_i)^2dy
=\frac{2n+2}{2n+4}\cdot\frac{(n!!)^2}{(n+1)!!(n-1)!!}\cdot\frac{1}{2n}\cdot\frac{\pi}{2}\varepsilon^2\equiv
M_o(n)\varepsilon^2.\nonumber
\end{equation}
\end{lem}

\begin{pf} The first three terms are obviously. Using left-invariance and symmetry, we have
\begin{align*}
\fint_{B_{\varepsilon}(x)}\psi(x^{-1}\circ
y)(y_i-x_i)^2dy=\frac{1}{2n}\fint_{B_{\varepsilon}(0)}\psi(y)\left|\overline{y}\right|^2dy.
\end{align*}
By using \eqref{polar2} and \eqref{jacobi}
\begin{align*}
\fint_{B_{\varepsilon}(0)}\psi(y)\left|\overline{y}\right|^2dy
=\frac{\int_{0}^{\varepsilon}\rho^{2n+3}d\rho\int_{0}^{\pi}\sin^{n+1}\phi
d\phi}
{\int_{0}^{\varepsilon}\rho^{2n+1}d\rho\int_{0}^{\pi}\sin^{n}\phi
d\phi}\\
=\frac{2n+2}{2n+4}\varepsilon^2\frac{\int_{0}^{\pi}\sin^{n+1}\phi
d\phi} {\int_{0}^{\pi}\sin^{n}\phi d\phi}.
\end{align*}
According to the integrals \begin{displaymath}
\int_{0}^{\frac{\pi}{2}}\sin^nxdx = \left\{
\begin{array}{lll}
            \frac{(2k-1)!!}{(2k)!!}\frac{\pi}{2},\hspace{2mm}n=2k,
           \\[2mm]
             \frac{(2k)!!}{(2k+1)!!}, \hspace{2mm}n=2k+1,
         \end{array}  \right.
\end{displaymath}
we obtain the desired results in this lemma.
\end{pf}

\section{Proof of Theorem \ref{thm2-1}} The general approach for the
proof of Theorem \ref{thm2-1} is similar to \cite{ES}; see also
\cite{JK}, \cite{manfredi2}. However, we notice that, the
Crandall-Ishii-Lions maximum principle (see Theorem 3.2 in
\cite{lions})  is not available for sub-parabolic structure on the
Heisenberg group. To circumvent this, one may use the Euclidean
Crandall-Ishii-Lions maximum principle to get the Euclidean jets,
and then twist the Euclidean jets to form sub-parabolic jets on
$\mathbb{H}^n$.  This method was introduced by Bieske \cite{bieske2}
for studying  existence and uniqueness of the viscosity solutions to
the sub-infinite Laplace equations on $\mathbb{H}^1$.
\begin{lem}(\cite{bieske2})\label{subjets}
Let $(a,p,Y)\in \mathbb{R}\times \mathbb{R}^{2n+1}\times
\mathcal{S}^{2n+1}$, and $\|\cdot\|_{E}$ denote the standard norm in
$\mathbb{R}^{2n+1}$. Define the standard Euclidean super-jet,
denoted by $\mathcal{J}^{2,+}_{E}$,
\begin{align}\label{euclidean-jet}
\mathcal{J}^{2,+}_{E}(u,(t^0,x^0))=\Big\{(a,p,Y): &u(t,x)\leq
u(t^0,x^0)+a(t-t^0)+\left<p,x-x^0\right>\nonumber\\
&+\frac{1}{2}\left<Y(x-x^0),(x-x^0)\right>+o\left(|t-t^0|+\|x-x^0\|_E^2\right)\Big\}.
\end{align}
denote
\begin{equation}\label{A}
A(x^0)=\left(\begin{array}{lll}\displaystyle 1\cdots 0 & 0 \cdots   0 & 2x^0_{n+1} \\
\displaystyle \vdots & \vdots & \vdots \\
\displaystyle 0 \cdots 1 & 0\cdots 0 & 2x^0_{2n}\\[2mm]
\displaystyle 0 \cdots 0 & 1 \cdots   0 & -2x^0_{1} \\
\displaystyle \vdots & \vdots & \vdots \\
\displaystyle 0 \cdots 0 & 0 \cdots   1& -2x^0_{n} \\
\displaystyle 0\cdots 0 & 0 \cdots   0 & 1 \\
\end{array}\right),
\end{equation}
Then
\begin{equation*}
(a,p,Y)\in \mathcal{J}^{2,+}_{E}\left(u,(t^0,x^0)\right)
\end{equation*}
implies
\begin{equation*}
\left(a,A(x^0)\cdot p,(A\cdot Y\cdot A^{T})_{2n}\right)\in
\mathcal{J}^{2,+}\left(u,(t^0,x^0)\right)
\end{equation*}
with the convention that for any matrix $M$, $M_{2n}$ is the
$2n\times2n$ principal minor.
\end{lem}

Now we are in a position to prove Theorem \ref{thm2-1}.

\textit{ { Proof of Theorem \ref{thm2-1}}}\quad Suppose $u$ is not a
viscosity super-solution of (\ref{operator1}) in the sense of
Definition \ref{def-p}, but satisfies the assumptions of  Theorem
\ref{thm2-1}. Then there exist $(t^0,x^0)\in \Omega_{T}$ and
$\phi\in C_{H}^2(\Omega_{T})$ satisfying $u(t^0,x^0)=\phi(t^0,x^0)$,
and $u(t,x)>\phi(t,x) \hspace{2mm}\mbox{for}
 \hspace{2mm}(t,x)\neq (t^0,x^0)$, for which $\nabla_0\phi(t^0,x^0)=0$,
 $(X^2\phi)^{*}(t^0,x^0)\neq0$, and
 \begin{equation}\label{ineq1}
 \phi_{t}(t^0,x^0)<
 \lambda_{min}\left((p-2)(X^2\phi)^{*}(t^0,x^0)\right)+\Delta_{H}\phi(t^0,x^0).
 \end{equation}
 Let
 \begin{equation}\label{w}
 w^{\alpha}(t,x,s,y)=u(t,x)-\phi(s,y)+\varphi(t,x,s,y),
\end{equation}
where
\begin{equation*}
 \varphi(t,x,s,y)=\frac{\alpha}{4}\left(\sum_{i=1}^{2n}\left(x_i-y_i\right)^4+\left(x_{2n+1}-y_{2n+1}
 +2\sum_{i=1}^{n}(x_{n+i}y_i-x_iy_{n+i})\right)^4\right)+\frac{\alpha}{2}(t-s)^2,
\end{equation*}
 and denote by $(t^{\alpha},x^{\alpha},s^{\alpha},y^{\alpha})$ the
 minimum point of $w^{\alpha}$ in $\overline{\Omega}_{T}\times\overline{\Omega}_{T}.$
 Since $(t^0,x^0)$ is a local minimum for $u-\phi$ and by
 \cite{bieske2}, we may assume that
 \begin{equation*}
 (t^{\alpha},x^{\alpha},s^{\alpha},y^{\alpha})\rightarrow
 (t^{0},x^{0},t^{0},x^{0})
 \hspace{3mm}\mbox{as}\hspace{1mm}\alpha\rightarrow +\infty.
 \end{equation*}
In particular, $(t^{\alpha},x^{\alpha})\in\Omega_{T}$  and
$(s^{\alpha},y^{\alpha})\in\Omega_{T}$ for all $\alpha$ large
enough.

We consider two cases: either
$\overline{x^{\alpha}}=\overline{y^{\alpha}}$  or
$\overline{x^{\alpha}}\neq\overline{y^{\alpha}}$ for all $\alpha$
large enough.

\textit{Case 1:} Let $\overline{x^{\alpha}}=\overline{y^{\alpha}}$,
and denote
\begin{equation}
\vartheta(s,y)=\varphi(t^{\alpha},x^{\alpha},s,y).
\end{equation}
Then $\phi(s,y)-\vartheta(s,y)$ has a local maximum at
$(s^{\alpha},y^{\alpha}),$ and  thus
\begin{equation*}
\phi_{s}(s^{\alpha},y^{\alpha})=\vartheta_{s}(s^{\alpha},y^{\alpha})
\hspace{1mm}\mbox{and}\hspace{1mm}
(X^2\phi)^{*}(s^{\alpha},y^{\alpha})\leq
(X^{2}\vartheta)^{*}(s^{\alpha},y^{\alpha}).
\end{equation*}
 A direct calculation yields
\begin{equation*}
(X^2\vartheta)^{*}(s^{\alpha},y^{\alpha})=0 \hspace{2mm}
\mbox{provided}\hspace{2mm}
\overline{x^{\alpha}}=\overline{y^{\alpha}}, \mbox{and}\hspace{2mm}
\vartheta_s(s^{\alpha},y^{\alpha})=-\alpha(t^{\alpha}-s^{\alpha}),
\end{equation*}
and thus,
\begin{equation}\label{maxpoint}
\phi_{s}(s^{\alpha},y^{\alpha})=-\alpha(t^{\alpha}-s^{\alpha}),
\hspace{1mm}\mbox{and}\hspace{1mm}
(X^2\phi)^{*}(s^{\alpha},y^{\alpha})\leq 0.
\end{equation}
By (\ref{ineq1}) and continuity of
\begin{equation*}
(t,x)\mapsto
\lambda_{min}\left((p-2)(X^2\phi)^{*}(t,x)\right)+\Delta_{H}\phi(t,x),
\end{equation*}
we have
\begin{equation}\label{ineq2}
\phi_{s}(s^{\alpha},y^{\alpha})<
\lambda_{min}\left((p-2)(X^2\phi)^{*}(s^{\alpha},y^{\alpha})\right)+\Delta_{H}\phi(s^{\alpha},y^{\alpha})
\end{equation}
for $\alpha$ large enough.
By (\ref{maxpoint}) and (\ref{ineq2}),
we have
\begin{equation}\label{ineq3}
0<-\vartheta_{s}(s^{\alpha},y^{\alpha})=\alpha(t^{\alpha}-s^{\alpha}),
\end{equation}
for $\alpha$ large enough. If $1<p<2$, the inequality follows from
the calculation
\begin{align}\label{same}
&\lambda_{min}\left((p-2)(X^2\phi)^{*}(s^{\alpha},y^{\alpha})\right)+\Delta_{H}\phi(s^{\alpha},y^{\alpha})\nonumber\\[1mm]
&=(p-2)\lambda_{max}\left((X^2\phi)^{*}(s^{\alpha},y^{\alpha})\right)
+trace\left((X^2\phi)^{*}(s^{\alpha},y^{\alpha})\right)\nonumber\\[1mm]
&=(p-1)\lambda_{max}+\sum_{\lambda_i\neq \lambda_{max}}\lambda_i\leq
0,
\end{align}
where $\lambda_i$, $\lambda_{max}$ denote the eigenvalue  and the
maximum eigenvalue  of $(X^2\phi)^{*}(s^{\alpha},y^{\alpha})$,
respectively.

  On the other hand, let \begin{equation}
  \mu(t,x)=-\varphi(t,x,s^{\alpha},y^{\alpha}).
 \end{equation}
 Similarly,
 $u(t,x)-\mu(t,x)$ has a local minimum at $(t^{\alpha},x^{\alpha})$,
 and
 \begin{equation*}
 \nabla_{0}\mu(t^{\alpha},x^{\alpha})=0,\hspace{1mm}
  (X^2\mu)^{*}(t^{\alpha},x^{\alpha})=0,
  \hspace{1mm}\mbox{provided}\hspace{1mm}
\overline{x^{\alpha}}=\overline{y^{\alpha}}.
\end{equation*}
 That is, $\mu$ is  a $C^2_{H}$ test function, by the
assumption on $u$, we have
\begin{equation}\label{ineq4}
0\leq \mu_{t}(t^{\alpha},x^{\alpha})=-\alpha(t^{\alpha}-s^{\alpha}),
\end{equation}
for $\alpha$ large enough. Summing up (\ref{ineq3}) and
(\ref{ineq4}) gives
\begin{equation*}
0<\alpha(t^{\alpha}-s^{\alpha})-\alpha(t^{\alpha}-s^{\alpha})=0.
\end{equation*}
This is a contradiction.

\textit{Case 2:} Next we consider the case
$\overline{x^{\alpha}}\neq\overline{y^{\alpha}}$ for all $\alpha$
large enough. We apply the Euclidean maximum   principle for
semi-continuous functions of Crandall-Ishii-Lions (see Theorem 3.2
in \cite{lions}).  There exists  $(2n+1)\times (2n+1)$ symmetric
matrices $Y^{\alpha}$, $Z^{\alpha}$ such that
\begin{align*}\displaystyle
\left(-D_s\varphi(t^{\alpha},x^{\alpha},s^{\alpha},y^{\alpha}),
-D_{y}\varphi(t^{\alpha},x^{\alpha},s^{\alpha},y^{\alpha}),
Y^{\alpha}\right)\in\overline{\mathcal{J}}^{2,+}_{E}\phi\left(s^{\alpha},y^{\alpha}\right),
\\
\displaystyle
\Big(D_t\varphi(t^{\alpha},x^{\alpha},s^{\alpha},y^{\alpha}),
D_{x}\varphi(t^{\alpha},x^{\alpha},s^{\alpha},y^{\alpha}),
Z^{\alpha}\Big)\in\overline{\mathcal{J}}_{E}^{2,-}u\left(t^{\alpha},x^{\alpha}\right).
\end{align*}
with the property that
\begin{equation}\label{property}
\left<Y^{\alpha}\gamma,\gamma\right>-\left<Z^{\alpha}\chi,\chi\right>\leq
\left<C\gamma\oplus\chi,\gamma\oplus\chi\right>,
\end{equation}
where
\begin{equation*}
C=B+\frac{1}{\alpha}B^2,\hspace{3mm} \gamma\oplus\chi=(\gamma,\chi),
\end{equation*}
and
\begin{equation*}
B=D^2_{x,y}\varphi\left(t^{\alpha},x^{\alpha},s^{\alpha},y^{\alpha}\right),
\end{equation*}
with the notations $D_x$, $D_y$ and $D_{xy}$ denote the Euclidean
derivatives.
 By using Lemma \ref{subjets} and the fact
 \begin{equation*}
 -D_s\varphi(t^{\alpha},x^{\alpha},s^{\alpha},y^{\alpha})=D_t\varphi(t^{\alpha},x^{\alpha},
 s^{\alpha},y^{\alpha})=\alpha(t^{\alpha}-s^{\alpha}),
 \end{equation*}
 we conclude that
\begin{eqnarray}\label{jet1}
\left(\alpha(t^{\alpha}-s^{\alpha}),-\nabla_{y}\varphi(t^{\alpha},x^{\alpha},
s^{\alpha},y^{\alpha}),\widetilde{Y^{\alpha}}\right)
\in\overline{\mathcal{J}}^{2,+}\phi(s^{\alpha},y^{\alpha}),
\end{eqnarray}
and
\begin{eqnarray}\label{jet2}
\left(\alpha(t^{\alpha}-s^{\alpha}),-\nabla_{x}\varphi(t^{\alpha},x^{\alpha},
s^{\alpha},y^{\alpha}),\widetilde{Z^{\alpha}}\right)
\in\overline{\mathcal{J}}^{2,-}u(t^{\alpha},x^{\alpha}),
\end{eqnarray}
where $\widetilde{Y^{\alpha}}$ and $\widetilde{Z^{\alpha}}$ are
$2n\times 2n$  symmetric matrices  defined by
\begin{eqnarray*}
\widetilde{Y^{\alpha}}=\left(A(y^{\alpha})\cdot Y^{\alpha}\cdot
A(y^{\alpha})^{T}\right)_{2n}
\end{eqnarray*}
and
\begin{eqnarray*}
\widetilde{Z^{\alpha}}=\left(A(x^{\alpha})\cdot Z^{\alpha}\cdot
A(x^{\alpha})^{T}\right)_{2n},
\end{eqnarray*}
where $A(x^{\alpha})$ and $A(y^{\alpha})$ are as in (\ref{A}) with
the point $x^{0}$ is replaced by $x^{\alpha}$ and $y^{\alpha}$,
respectively.

\textit{ Claim:}  Let $\xi=\overline{(y^{\alpha})^{-1}\circ
x^{\alpha}}\in \mathbb{R}^{2n}$, we have the following estimate
\begin{eqnarray}\label{claim}
\left<\widetilde{Y^{\alpha}}\xi,\xi\right>-\left<\widetilde{Z^{\alpha}}\xi,\xi\right>\leq
0, \hspace{2mm} \mbox{as} \hspace{2mm}\alpha\rightarrow +\infty.
\end{eqnarray}

Assume the above claim is true. By (\ref{ineq1}), there exists a
constant $\theta>0$, such that
 \begin{equation}
\theta+\phi_{t}(t^0,x^0)<\lambda_{min}\left((p-2)(X^2\phi)^{*}(t^0,x^0)\right)+\Delta_{H}\phi(t^0,x^0),
 \end{equation}
 and with the continuity of
\begin{equation*}
(t,x)\mapsto
\lambda_{min}\left((p-2)(X^2\phi)^{*}(t,x)\right)+\Delta_{H}\phi(t,x),
\end{equation*}
we have
\begin{equation}\label{ineq-theta}
\theta+\phi_{s}(s^{\alpha},y^{\alpha})<
\lambda_{min}\left((p-2)(X^2\phi)^{*}(s^{\alpha},y^{\alpha})\right)+\Delta_{H}\phi(s^{\alpha},y^{\alpha})
\end{equation}
for $\alpha$ large enough.

 Using (\ref{jet1}), (\ref{jet2}),
(\ref{ineq-theta}) and the the assumptions on $u$, we have
\begin{align*}
\theta&=\theta+\alpha(t^{\alpha}-s^{\alpha})-\alpha(t^{\alpha}-s^{\alpha})\\[2mm]
&<
(p-2)\left<\widetilde{Y^{\alpha}}\frac{\overline{(y^{\alpha})^{-1}\circ
x^{\alpha}}}{|(y^{\alpha})^{-1}\circ
x^{\alpha}|},\frac{\overline{(y^{\alpha})^{-1}\circ
x^{\alpha}}}{|(y^{\alpha})^{-1}\circ
x^{\alpha}|}\right>+trace(\widetilde{Y^{\alpha}})\\[2mm]
&-(p-2)\left<\widetilde{Z^{\alpha}}\frac{\overline{(y^{\alpha})^{-1}\circ
x^{\alpha}}}{|(y^{\alpha})^{-1}\circ
x^{\alpha}|},\frac{\overline{(y^{\alpha})^{-1}\circ
x^{\alpha}}}{|(y^{\alpha})^{-1}\circ
x^{\alpha}|}\right>-trace(\widetilde{Z^{\alpha}})\\[2mm]
&=(p-2)\left<(\widetilde{Y^{\alpha}}-\widetilde{Z^{\alpha}})\frac{\overline{(y^{\alpha})^{-1}\circ
x^{\alpha}}}{|(y^{\alpha})^{-1}\circ
x^{\alpha}|},\frac{\overline{(y^{\alpha})^{-1}\circ
x^{\alpha}}}{|(y^{\alpha})^{-1}\circ
x^{\alpha}|}\right>+trace(\widetilde{Y^{\alpha}}-\widetilde{Z^{\alpha}})\\[2mm]
&\leq 0,
\end{align*}
the last inequality being valid by the claim (\ref{claim}) in the
case $p>2$. If $1<p<2$, the last inequality follows the same
calculation in (\ref{same}).
This is a contradiction.

\textit {Proof of the claim.} For any $\xi=(\xi_1,\cdots,\xi_{2n})\in\mathbb{R}^{2n}$, let
\begin{equation}\label{gamma}
\zeta=\left(\xi,2\sum_{i=1}^{n}(\xi_iy^{\alpha}_{n+i}-\xi_{n+i}y^{\alpha}_{i})\right),\quad
\eta=\left(\xi,2\sum_{i=1}^{n}(\xi_ix^{\alpha}_{n+i}-\xi_{n+i}x^{\alpha}_{i})\right).
\end{equation}
Recalling the definitions of $\widetilde{Y^{\alpha}}$,
$\widetilde{Z^{\alpha}}$, and combining (\ref{property}), we obtain
\begin{equation*}
\left<\widetilde{Y^{\alpha}}\xi,\xi\right>-\left<\widetilde{Z^{\alpha}}\xi,\xi\right>
=\left<Y^{\alpha}\zeta,\zeta\right>-\left<Z^{\alpha}\eta,\eta\right>
\leq\left<C\zeta\oplus\eta,\zeta\oplus\eta\right>.
\end{equation*}
Straightforward computations show that
\begin{equation}
\left<B\zeta\oplus\eta,\zeta\oplus\eta\right>=0,
\end{equation}
and
\begin{equation}
\left<B^2\zeta\oplus\eta,\zeta\oplus\eta\right>=8\alpha^2\xi^2\left(x_{2n+1}^{\alpha}-y_{2n+1}^{\alpha}
 +2\sum_{i=1}^{n}(x_{n+i}^{\alpha}y_i^{\alpha}-x_i^{\alpha}y_{n+i}^{\alpha})\right)^6.
\end{equation}
Now choosing $\xi=\overline{(y^{\alpha})^{-1}\circ
x^{\alpha}}=(x_1^{\alpha}-y_1^{\alpha},\cdots,x_{2n}^{\alpha}-y_{2n}^{\alpha})$, and noting that \cite{Bieske2006}
\begin{equation*}
\lim_{\alpha\rightarrow+\infty}\alpha\left(\sum_{i=1}^{2n}\left(x_i-y_i\right)^4+\left(x_{2n+1}-y_{2n+1}
 +2\sum_{i=1}^{n}(x_{n+i}y_i-x_iy_{n+i})\right)^4\right)=0
\end{equation*}
Thanks to  $C=B+1/\alpha B^2$,  we have
\begin{equation}
\left<C\zeta\oplus\eta,\zeta\oplus\eta\right>\rightarrow 0,
\hspace{2mm} \mbox{as}\hspace{2mm}\alpha\rightarrow +\infty.
\end{equation}
The claimed (\ref{claim}) is proved.
\qed \\

\section{Proof of the main results}
In this section, we  prove  asymptotic mean value formulae for
sub-heat equations (i.e. $p=2$) and sub-infinity Laplace parabolic
equations (i.e. $p=\infty$) on $\mathbb{H}^n$,   and construct an
example to show that the formulae do not hold in non-asymptotic
sense. We begin with a key lemma, which depicts the directions of
horizontal maximum and minima of a function, and whose Euclidean
version is obvious (cf. \cite{LW08}).

For $\phi\in C^2_{H}(\Omega),$ $x^0\in\Omega$ and $r>0$ with
$\overline{B_{r}(x^0)}\subset\Omega$, we define
\begin{equation*}
M(r)=\max\limits_{x\in\partial B_r(x^0)}\phi(x), \hspace{3mm}
\mbox{and} \hspace{3mm} m(r)=\min\limits_{x\in\partial
B_r(x^0)}\phi(x).
\end{equation*}
In addition, $(x^{r})^{+}\in\partial B_r(x^0)$ and
$(x^{r})^{-}\in\partial B_r(x^0)$ denote any point  such that
\begin{equation*}
\phi\left((x^{r})^+\right)=M(r),\hspace{3mm} \mbox{and}\hspace{3mm}
\phi\left((x^{r})^-\right)=m(r).
\end{equation*}
 Define the set of horizontal
maximum  directions of $\phi$ at $x^0$ to be the set
\begin{displaymath}
E^{+}(x^0)=\Big\{\lim_{k}\frac{\overline{(x^0)^{-1}\circ
(x^{r_k})^{+}}}{r_k} \mbox{ \hspace{1mm}for some  sequence
}\hspace{1mm} r_k\rightarrow 0\Big\},
\end{displaymath}
and the set of horizontal minimum directions of $\phi$ at $x^0$ to
be the set
\begin{displaymath}
E^{-}(x^0)=\Big\{\lim_{k}\frac{\overline{(x^0)^{-1}\circ
(x^{r_k})^{-}}}{r_k} \mbox{ \hspace{1mm}for some  sequence
}\hspace{1mm} r_k\rightarrow 0\Big\}.
\end{displaymath}

\begin{lem}\label{lemma3-1}
Let $\phi\in C_{H}^2$ and $\nabla_0\phi(x^0)\neq 0$, then
$$E^{+}(x^0)=\frac{\nabla_0\phi}{|\nabla_0\phi|}(x^0),\hspace{3mm}\mbox{and}
\hspace{3mm}E^{-}(x^0)=-\frac{\nabla_0\phi}{|\nabla_0\phi|}(x^0).$$
\end{lem}

\begin{pf} Define a Lagrange function to be
\begin{align}
F(x)&=\phi(x)+\lambda\left(\rho^4((x^0)^{-1}\circ
x)-\varepsilon^4\right)\nonumber\\
&=\phi(\overline{x},x_{2n+1})+\lambda\Big\{|\overline{x}-\overline{x^{0}}|^4
+\Big(x_{2n+1}-x^0_{2n+1}+2\sum_{i=1}^n(x^0_ix_{n+i}-x_ix^0_{n+i})\Big)^2-\varepsilon^4\Big\}.\nonumber
\end{align}
If $x^{\varepsilon}$ is a solution of  $\min\limits_{\partial
B_{\varepsilon}(x^0)}\phi(x)$, then there exists
$\lambda^{\varepsilon}$, such that for $i=1,\cdots,n$
\begin{displaymath}
\left\{
\begin{array}{lll}
0&=X_{i}F(x^{\varepsilon})\\
&=X_{i}\phi(x^{\varepsilon})
+4\lambda^{\varepsilon}\Big\{\left|\overline{x^{\varepsilon}}-\overline{x^0}\right|^2\left(x^{\varepsilon}_i-x^0_i\right)
+\left(x^{\varepsilon}_{2n+1}-x^0_{2n+1}+2\sum_{i=1}^n\left(x^0_ix^{\varepsilon}_{n+i}-x^{\varepsilon}_ix^0_{n+i}\right)\right)\cdot
\left(x^{\varepsilon}_{n+i}-x^0_{n+i}\right)\Big\},\\[3mm]
0&=X_{n+i}F(x^{\varepsilon})\\
&=X_{n+i}\phi(x^{\varepsilon})
+4\lambda^{\varepsilon}\Big\{\left|\overline{x^{\varepsilon}}-\overline{x^0}\right|^2\left(x^{\varepsilon}_{n+i}-x^0_{n+i}\right)
+\left(x^{\varepsilon}_{2n+1}-x^0_{2n+1}+2\sum_{i=1}^n\left(x^0_ix^{\varepsilon}_{n+i}-x^{\varepsilon}_ix^0_{n+i}\right)\right)\cdot
\left(x^0_{i}-x^{\varepsilon}_{i}\right)\Big\},\\[3mm]
0&= TF(x^{\varepsilon})\\
&=T\phi(x^{\varepsilon})
+2\lambda^{\varepsilon}\left(x^{\varepsilon}_{2n+1}-x^0_{2n+1}
+2\sum_{i=1}^n\left(x^0_ix^{\varepsilon}_{n+i}-x^{\varepsilon}_ix^0_{n+i}\right)\right),\\[3mm]
0&=\frac{\partial F}{\partial
\lambda}(x^{\varepsilon})\\
&=\left|\overline{x^{\varepsilon}}-\overline{x^0}\right|^4+\left(x^{\varepsilon}_{2n+1}-x^0_{2n+1}
+2\sum_{i=1}^n\left(x^0_ix^{\varepsilon}_{n+i}-x^{\varepsilon}_ix^0_{n+i}\right)\right)^2-\varepsilon^4.
\end{array}\right.
\end{displaymath}
A direct computation yields
\begin{displaymath}
 \Big|\nabla_0\phi\Big|(x^{\varepsilon})=\sum_{i=1}^{2n}X_i\phi\left(x^{\varepsilon}\right)
 =4\lambda^{\varepsilon}\varepsilon^2\left|\overline{x^{\varepsilon}}-\overline{x^0}\right|.
 \end{displaymath}
Therefore
\begin{displaymath}\label{X1}
\frac{X_i\phi}{|\nabla_0\phi|}\Big(x^{\varepsilon}\Big)
=-\frac{\left|\overline{x^{\varepsilon}}-\overline{x^0}\right|^2\left(x^{\varepsilon}_i-x^0_i\right)
+\left(x^{\varepsilon}_{2n+1}-x^0_{2n+1}+2\sum_{i=1}^n\left(x^0_ix^{\varepsilon}_{n+i}-x^{\varepsilon}_ix^0_{n+i}\right)\right)
\cdot\left(x^{\varepsilon}_{n+i}-x^0_{n+i}\right)}
{\varepsilon^2\left|\overline{x^{\varepsilon}}-\overline{x^0}\right|}.
 \end{displaymath}
Similarly, \begin{eqnarray}\label{X2}
\frac{X_{n+i}\phi}{|\nabla_0\phi|}\Big(x^{\varepsilon}\Big) =-
\frac{\left|\overline{x^{\varepsilon}}-\overline{x^0}\right|^2\left(x^{\varepsilon}_{n+i}-x^0_{n+i}\right)
+\left(x^{\varepsilon}_{2n+1}-x^0_{2n+1}+2\sum_{i=1}^n\left(x^0_ix^{\varepsilon}_{n+i}-x^{\varepsilon}_ix^0_{n+i}\right)\right)
\cdot\left(x^0_{i}-x^{\varepsilon}_{i}\right)}
{\varepsilon^2\left|\overline{x^{\varepsilon}}-\overline{x^0}\right|},\nonumber
 \end{eqnarray}
and
\begin{equation}\label{T}
\frac{T\phi}{|\nabla_0\phi|}\Big(x^{\varepsilon}\Big)
=-\frac{x^{\varepsilon}_{2n+1}-x^0_{2n+1}+2\sum_{i=1}^n\left(x^0_ix^{\varepsilon}_{n+i}-x^{\varepsilon}_ix^0_{n+i}\right)}
{2\varepsilon^2\left|\overline{x^{\varepsilon}}-\overline{x^0}\right|}.
 \end{equation}
Let $\varepsilon\rightarrow 0$ in (\ref{T}),  we get
\begin{equation}\label{t}
\frac{T\phi}{|\nabla_0\phi|}\left(x^0\right)
=\lim_{\varepsilon\rightarrow0}-\frac{x^{\varepsilon}_{2n+1}-x^0_{2n+1}+2\sum_{i=1}^n
\left(x^0_ix^{\varepsilon}_{n+i}-x^{\varepsilon}_ix^0_{n+i}\right)}
{2\varepsilon^2\left|\overline{x^{\varepsilon}}-\overline{x^0}\right|}.\nonumber
\end{equation}
Therefore
\begin{equation}
\lim_{\varepsilon\rightarrow0}\frac{\left(x^{\varepsilon}_{2n+1}-x^0_{2n+1}+2\sum_{i=1}^n\left(x^0_ix^{\varepsilon}_{n+i}
-x^{\varepsilon}_ix^0_{n+i}\right)\right)\left(x^{\varepsilon}_{n+i}-x^0_{n+i}\right)}
{\varepsilon^2\left|\overline{x^{\varepsilon}}-\overline{x^0}\right|}=0,\nonumber
\end{equation}

\begin{equation}
\lim_{\varepsilon\rightarrow0}\frac{\left(x^{\varepsilon}_{2n+1}-x^0_{2n+1}+2\sum_{i=1}^n\left(x^0_ix^{\varepsilon}_{n+i}
-x^{\varepsilon}_ix^0_{n+i}\right)\right)\left(x^0_{i}-x^{\varepsilon}_{i}\right)}
{\varepsilon^2|\overline{x^{\varepsilon}}-\overline{x^0}|}=0,\nonumber
\end{equation}
and
\begin{displaymath}
\lim_{\varepsilon\rightarrow0}\frac{|\overline{x^{\varepsilon}}-\overline{x^0}|}{\varepsilon}=1.
\end{displaymath}
Hence
\begin{displaymath}
-\frac{X_i\phi}{|X\phi|}(x^0)
=\lim_{\varepsilon\rightarrow0}\frac{x^{\varepsilon}_i-x^0_i}{\varepsilon}
\cdot\lim_{\varepsilon\rightarrow0}\frac{|\overline{x^{\varepsilon}}-\overline{x^0}|}
{\varepsilon}=\lim_{\varepsilon\rightarrow0}\frac{x^{\varepsilon}_i-x^0_i}{\varepsilon},
 \end{displaymath}
and
 \begin{displaymath}
-\frac{X_{n+i}\phi}{|X\phi|}\Big(x^0\Big)
=\lim_{\varepsilon\rightarrow0}\frac{x^{\varepsilon}_{n+i}-x^0_{n+i}}{\varepsilon}
\cdot\lim_{\varepsilon\rightarrow0}\frac{|\overline{x^{\varepsilon}}-\overline{x^0}|}
{\varepsilon}=\lim_{\varepsilon\rightarrow0}\frac{x^{\varepsilon}_{n+i}-x^0_{n+i}}{\varepsilon}.
\end{displaymath}
That is $E^{-}(x^0)=-\frac{X\phi}{|\nabla_0\phi|}(x^0)$. The same
argument to show $E^{+}(x^0)=\frac{X\phi}{|\nabla_0\phi|}(x^0).$
Therefore, the proof of the lemma is complete.
\end{pf}
\vspace{3mm}

Now, we prove an asymptotic mean value formula of the sub-heat
equations on $\mathbb{H}^{n}$.

 \textit{Proof of Theorem \ref{thm-p=2}:}
Let $u$ be a smooth function, and $(t,x)\in\Omega_{T}$. Consider the
Taylor expansion
\begin{align}\label{taylor}
u(s,y)&=u(t,x)+u_t(t,x)(s-t)+\nabla
u(t,x)\cdot(x^{-1}\circ y)\nonumber\\
&+\frac{1}{2}\left<(X^2 u)^{*}(t,x)\overline{(x^{-1}\circ
y)},\overline{(x^{-1}\circ y)}\right> +c\left(\rho^2(x^{-1}\circ
y)+|s-t|\right).
\end{align}
Averaging both sides of (\ref{taylor}), we have
\begin{align}\label{equ1}
&\fint_{t-\varepsilon^2}^{t}\fint_{B_{\varepsilon}(x)}\psi(x^{-1}\circ
y)u(s,y)dyds\nonumber\\
&=u(t,x)+u_t(t,x)\fint_{t-\varepsilon^2}^{t}(s-t)d s
+\fint_{B_{\varepsilon}(x)}\psi(x^{-1}\circ y)\nabla
u(t,x)\cdot(x^{-1}\circ y)dy \nonumber\\
&+\frac{1}{2}\fint_{B_{\varepsilon}(x)}\psi(x^{-1}\circ y)\left<(X^2
u)^{*}(t,x)\overline{(x^{-1}\circ y)},\overline{(x^{-1}\circ
y)}\right>dy + o(\varepsilon^2).
\end{align}
By Lemma \ref{lem2}, we get
\begin{equation}\label{3-4}
\fint_{B_{\varepsilon}(x)}\psi(x^{-1}\circ y)\nabla
u(t,x)\cdot(x^{-1}\circ y)dy=0,
\end{equation}
and
\begin{equation}\label{4}
\frac{1}{2}\fint_{B_{\varepsilon}(x)}\psi(x^{-1}\circ
y)\left<(X^2u)^{*}(t,x)\overline{(x^{-1}\circ
y)},\overline{(x^{-1}\circ y)}\right>dy
=\frac{1}{2}M(n)\varepsilon^2\Delta_{H}u(t,x).
\end{equation}
Finally,
\begin{equation}\label{5}
\fint_{t-\varepsilon^2}^{t}(s-t)ds=-\frac{1}{2}\varepsilon^2.
\end{equation}
Substituting (\ref{3-4}), (\ref{4}) and (\ref{5}) into (\ref{equ1}),
we have
\begin{equation}\label{identity}
\fint_{t-\varepsilon^2}^{t}\fint_{B_{\varepsilon}(x)}\psi(x^{-1}\circ
y)u(s,y)dyds=u(t,x)+\frac{1}{2}\varepsilon^2\Big(M(n)\Delta_{H}u(t,x)-u_t(t,x)\Big)+o(\varepsilon^2).
\end{equation}
This holds for any smooth function.

We first prove that if  $u$ satisfies the asymptotic mean value
formula (\ref{vis-p=2}), then $u$ is a  solution to
(\ref{equationp=2}). By (\ref{identity}), we have
\begin{align*}
u(t,x)&=\fint_{t-\varepsilon^2}^{t}\fint_{B_{\varepsilon}(x)}\psi(x^{-1}\circ
y)u(s,y)dyds+o(\varepsilon^2)\\
&=u(t,x)+\frac{1}{2}\varepsilon^2\Big(M(n)\Delta_{H}u(t,x)-u_t(t,x)\Big)+o(\varepsilon^2).
\end{align*}
That is
\begin{equation}\label{equ4}
\frac{1}{2}\varepsilon^2\Big(M(n)\Delta_{H}u(t,x)-u_t(t,x)\Big)+o(\varepsilon^2)=0.
\end{equation}
Dividing (\ref{equ4}) by $\varepsilon^2$ and passing to the limit
$\varepsilon\rightarrow 0$, we have
\begin{equation}
u_t(t,x)=M(n)\Delta_{H}u(t,x).
\end{equation}

Next we are ready to prove the converse implication. If  $u$ is a
solution of (\ref{equationp=2}), then (\ref{identity}) implies that
\begin{equation*}
u(t,x)=\fint_{t-\varepsilon^2}^{t}\fint_{B_{\varepsilon}(x)}\psi(x^{-1}\circ
y)u(s,y)dyds+o(\varepsilon^2).
\end{equation*}
This ends the proof.
 \qed

\vspace{3mm}
 The same argument shows that solutions to  the sub-heat
equation
\begin{equation*}
u_t(t,x)=\Delta_{H}u(t,x),
\end{equation*}
are characterized by the asymptotic mean value formula
\begin{equation}\label{vis2}
u(t,x)=\fint_{t-M(n)\varepsilon^2}^{t}
\fint_{B_{\varepsilon}(x)}\psi(x^{-1}\circ
y)u(s,y)dyds+o(\varepsilon^2)\hspace{2mm} \mbox{as}\hspace{2mm}
\varepsilon\rightarrow 0.
\end{equation}

Consider the mean value formula (\ref{p=2}) for $H-$harmonic
functions on $\mathbb{H}^n$,
 it is natural to ask if the formula
(\ref{vis2}) holds in a non-asymptotic sense. To be more precise, if
$u$ is  a  solution to
\begin{equation*}
u_t(t,x)=\Delta_{H}u(t,x),
\end{equation*}
does the equation
\begin{equation*}
u(t,x)=\fint_{t-M(n)\varepsilon^2}^{t}
\fint_{B_{\varepsilon}(x)}\psi(x^{-1}\circ y)u(s,y)dyds
\end{equation*}
 hold at all $(t,x)\in\Omega_{T}$ for all $\varepsilon>0$ enough
small. The answer to this question is negative, we give an example
as follows.

Let
$$u(t,x)=12t^2+12x_1^2t+x_1^4,$$
where $x=(x_1,x_2,x_3)\in \mathbb{H}^1.$ It is easy to check that
$u$ is a solution of
$$u_t(t,x)=\Delta_{H}u(t,x).$$
A direct calculation yields $M(1)=\frac{\pi}{12}$, and
\begin{equation*}
\int_{B_{\varepsilon}(0)}\psi(y)dy=\pi\varepsilon^4.
\end{equation*}
Thus
\begin{align*}
\fint_{B_{\varepsilon}(0)}\psi(y)u(s,y)dy&=\fint_{B_{\varepsilon}(0)}\psi(y)(12s^2+12y_1^2s+y_1^4)dy\\
&=12s^2+\pi\varepsilon^2s+\frac{1}{8}\varepsilon^4,
\end{align*}
and
\begin{align*}
&\fint_{1-\frac{\pi}{12}\varepsilon^2}^{1}(12s^2+\pi\varepsilon^2s+\frac{1}{8}\varepsilon^4)ds\\
&=12-\pi\varepsilon^2+\frac{1}{8}\varepsilon^4+\pi\varepsilon^4+\frac{1}{36}\pi^2\varepsilon^4-\frac{1}{24}\pi^2\varepsilon^6.
\end{align*}
 That is
\begin{equation*}
\fint_{1-\frac{\pi}{12}\varepsilon^2}^{1}
\fint_{B_{\varepsilon}(0)}\psi(y)u(y,s)dyds\neq u(0,1)=12.
\end{equation*}

\vspace{3mm}

Next, we characterize the viscosity solutions  of the homogeneous
sub-infinity Laplace parabolic equation in terms of an asymptotic
mean value formula on $\mathbb{H}^n$.

{ \textit{Proof of Theorem \ref{thm-infinity}}}\quad Choose a point
$(t,x)\in \Omega_{T}$, $\varepsilon>0$, $s\in (t-\varepsilon^2,t)$
and any $\phi\in C_{H}^2(\Omega_T)$. Denote by $x^{\varepsilon,s}$
be a point at which $\phi$ attains its minimum in
$\overline{B}_{\varepsilon}(x)$ at time $s$, that is
\begin{equation*}
\phi(s,x^{\varepsilon,s})=\min_{y\in\overline{B}_{\varepsilon}(x)}\phi(s,y).
\end{equation*}
Consider the Taylor expansion
\begin{align}\label{taylor2}
\phi(s,y)&=\phi(t,x)+\phi_t(t,x)(s-t)+\nabla
\phi(t,x)\cdot(x^{-1}\circ y)\nonumber\\
&+\frac{1}{2}\left<(X^2 \phi)^{*}(t,x)\overline{(x^{-1}\circ
y)},\overline{(x^{-1}\circ y)}\right> +c\left(\rho^2(x^{-1}\circ
y)+|s-t|\right).
\end{align}
 Taking $y=x^{\varepsilon,s}$ in(\ref{taylor2}) and noting
\begin{equation*}
x^{-1}\circ
x^{\varepsilon,s}=\left(x_1^{\varepsilon,s}-x_1,\cdots,x_{2n}^{\varepsilon,s}-x_{2n},x_{2n+1}^{\varepsilon,s}-x_{2n+1}
+2\sum\limits_{i=1}^{n}\left(x_{n+i}^{\varepsilon,s}x_i-x_{i}^{\varepsilon,s}x_{n+i}\right)\right),
\end{equation*}
 we have
\begin{align}\label{6}
\phi(s,x^{\varepsilon,s})&=\phi(t,x)+\phi_t(t,x)(s-t)+\nabla
\phi(t,x)(x^{-1}\circ x^{\varepsilon,s})
\nonumber\\
&+\frac{1}{2}\left<(X^2\phi)^{*}(t,x)\overline{(x^{-1}\circ
x^{\varepsilon,s})},\overline{(x^{-1}\circ x^{\varepsilon,s})}\right>+c\left(\varepsilon^2+|s-t|)\right)\nonumber\\
&=\phi(t,x)+\phi_t(t,x)(s-t)+\sum\limits_{i=1}^{2n}X_i\phi(t,x)(x_i^{\varepsilon,s}-x_i)\nonumber\\
&+T\phi(t,x)\left(x_{2n+1}^{\varepsilon,s}-x_{2n+1}
+2\sum\limits_{i=1}^{n}(x_{n+i}^{\varepsilon,s}x_i-x_{i}^{\varepsilon,s}x_{n+i})\right)\nonumber\\
&+\frac{1}{2}\sum_{i,j=1}^{2n}X_iX_j\phi(t,x)\cdot(x_{i}^{\varepsilon,s}-x_i)\cdot(x_{j}^{\varepsilon,s}-x_j)
+c\left(\varepsilon^2+|s-t|\right)\quad \hbox{as}\hspace{1mm}
\varepsilon\rightarrow 0.
\end{align}  Similarly, taking
$y=y^{\varepsilon,s}=\left(2x_1-x_1^{\varepsilon,s},\cdots,2x_{2n}-x_{2n}^{\varepsilon,s},2x_{2n+1}-x_{2n+1}^{\varepsilon,s}\right)$
in (\ref{taylor2}), and
\begin{equation*}
(x)^{-1}\circ
y^{\varepsilon,s}=\left(x_1-x_1^{\varepsilon,s},\cdots,x_{2n}-x_{2n}^{\varepsilon,s},
x_{2n+1}-x_{2n+1}^{\varepsilon,s}+2\sum\limits_{i=1}^{n}(x_{i}^{\varepsilon,s}x_{n+i}-x_ix_{n+i}^{\varepsilon,s})\right),
\end{equation*}
we have
\begin{align}\label{7}
\phi(s,y^{\varepsilon,s})
&=\phi(t,x)+\phi_t(t,x)(s-t)-\sum_{i=1}^{2n}X_i\phi(t,x)(x_i^{\varepsilon,s}-x_i)\nonumber\\
&-T\phi(t,x)\left(x_{2n+1}^{\varepsilon,s}-x_{2n+1}
+2\sum\limits_{i=1}^{n}(x_{n+i}^{\varepsilon,s}x_i-x_{i}^{\varepsilon,s}x_{n+i})\right)\nonumber\\
&+\frac{1}{2}\sum_{i,j=1}^{2n}X_iX_j\phi(t,x)\cdot(x_{i}^{\varepsilon,s}-x_i)\cdot(x_{j}^{\varepsilon,s}-x_j)
+c\left(\varepsilon^2+|s-t|\right).
\end{align}
Summing (\ref{6}) and (\ref{7}), we have
\begin{align*}
&\phi(s,x^{\varepsilon,s})+\phi(s,y^{\varepsilon,s})-2\phi(t,x)\\
&=2\phi_t(t,x)(s-t)+\sum_{i,j=1}^{2n}X_iX_j\phi(t,x)\cdot(x_{i}^{\varepsilon,s}-x_i)\cdot(x_{j}^{\varepsilon,s}-x_j)
+o\left(\varepsilon^2+|s-t|\right).
\end{align*}
Since $x^{\varepsilon,s}$ is a  minimum point of $\phi(\cdot,s)$ on
$\overline{B}_{\varepsilon}(x)$, we get
\begin{equation*}
\phi(s,x^{\varepsilon,s})+\phi(s,y^{\varepsilon,s})-2\phi(t,x)
\leq\max_{y\in\overline{B}_{\varepsilon}(x)}u(s,y)
+\min_{y\in\overline{B}_{\varepsilon}(x)}u(s,y)-2\phi(t,x),
\end{equation*}
and thus
\begin{align*}
&\max_{y\in\overline{B}_{\varepsilon}(x)}u(s,y)
+\min_{y\in\overline{B}_{\varepsilon}(x)}u(s,y)-2\phi(t,x)\\
&\geq
2\phi_t(x,t)(s-t)+\sum_{i,j=1}^{2n}X_iX_j\phi(t,x)\cdot(x_{i}^{\varepsilon,s}-x_i)\cdot(x_{j}^{\varepsilon,s}-x_j)
+o\left(\varepsilon^2+|s-t|\right).
\end{align*}
Integration over the time interval and the fact
$\fint_{t-\varepsilon^2}^{t}(s-t)ds=-\frac{1}{2}\varepsilon^2$ imply
\begin{align}\label{ineq6}
&\frac{1}{2}\fint_{t-\varepsilon^2}^{t}\left(\max_{y\in\overline{B}_{\varepsilon}(x)}\phi(s,y)
+\min_{y\in\overline{B}_{\varepsilon}(x)}\phi(y,s)\right)ds-\phi(t,x)\nonumber\\
&\geq\frac{\varepsilon^2}{2}\left(\fint_{t-\varepsilon^2}^{t}\sum_{i,j=1}^{2n}X_iX_j\phi(t,x)\cdot\frac{(x_{i}^{\varepsilon,s}-x_i)}{\varepsilon}
\cdot\frac{(x_{j}^{\varepsilon,s}-x_j)}{\varepsilon}ds
-\phi_t(t,x)\right)+o(\varepsilon^2).
\end{align}
This inequality holds for any function $\phi\in
C_{H}^2(\Omega_{T}).$

In the following, we prove the result via a dichotomy.

 Because $\phi\in C_{H}^2(\Omega),$ if
$\nabla_{0}\phi(t,x)\neq 0$, so $\nabla_{0}\phi(s,x)\neq 0$ for
$t-\varepsilon^2\leq s\leq t$ and for small enough $\varepsilon>0$,
and thus $x^{\varepsilon,s}\in\partial B_{\varepsilon}(x)$ for small
$\varepsilon$. By Lemma \ref{lemma3-1}, we have
\begin{equation}
\lim_{\varepsilon\rightarrow
0}\frac{x_i^{\varepsilon,s}-x_i}{\varepsilon}=-\frac{X_i\phi}{|\nabla_{0}\phi|}\Big(t,x\Big)\hspace{3mm}
\mbox{for} \hspace{1mm}i=1,\cdots,2n.
\end{equation}
Hence, we get the limit
\begin{align}\label{equ2}
\lim_{\varepsilon\rightarrow
0}\fint_{t-\varepsilon^2}^{t}\sum_{i,j=1}^{2n}X_iX_j\phi(t,x)\cdot\frac{(x_{i}^{\varepsilon,s}-x_i)}{\varepsilon}
\cdot\frac{(x_{j}^{\varepsilon,s}-x_j)}{\varepsilon}ds\nonumber\\
=\sum_{i,j=1}^{2n}X_iX_j\phi(t,x)
\cdot\frac{X_i\phi}{|\nabla_{0}\phi|}(t,x)\cdot\frac{X_j\phi}{|\nabla_{0}\phi|}(t,x)=\Delta_{H}^{\infty}\phi(t,x).
\end{align}

We first prove that if the asymptotic mean value formula
(\ref{vis-infinity}) holds for $u$ in viscosity sense, then $u$
satisfies the definition of viscosity solutions to
(\ref{parabolic-infinity}) whenever $\nabla_0\phi\neq0$. Let
$\phi\in C_{H}^2(\Omega_{T})$ be a test function such that $u-\phi$
has a strict minimum at $(t^0,x^0)$ and $\nabla_{0}\phi(t^0,x^0)\neq
0$, we have
\begin{equation}\label{ineq7}
0\geq-\phi(t^0,x^0)+\frac{1}{2}\fint_{t^0-\varepsilon^2}^{t^0}\left(\max_{y\in\overline{B_{\varepsilon}}(x^0)}\phi(s,y)
+\min_{y\in\overline{B_{\varepsilon}}(x^0)}\phi(s,y)\right)ds+o(\varepsilon^2).
\end{equation}
By (\ref{ineq6}), (\ref{equ2}) and (\ref{ineq7}), we have
\begin{equation}
o(\varepsilon^2)\geq
\frac{\varepsilon^2}{2}\left(\Delta_{H}^{\infty}\phi(t^0,x^0)-\phi_t(t^0,x^0)\right)+o(\varepsilon^2).
\end{equation}
Dividing by $\varepsilon^2$ and passing to  a limit, we get
\begin{equation}
\phi_{t}(t^0,x^0)\geq \Delta_{H}^{\infty}\phi(t^0,x^0).
\end{equation}
That is $u$ is a viscosity super-solution of
(\ref{parabolic-infinity}).

 To prove that $u$ is a viscosity sub-solution, we
first derive a reverse inequality to (\ref{ineq6}) by considering
the maximum point of $\phi$, and then we choose a test function
$\phi$ that touches $u$ from above.

To prove the reverse implication, assume that $u$ is a viscosity
super-solution of (\ref{parabolic-infinity}). Let $\phi\in
C_{H}^2(\Omega_{T})$ be a test function such that $u-\phi$ has a
strict minimum at $(t^0,x^0)$ and $\nabla_{0}\phi(t^0,x^0)\neq 0$,
we have
\begin{equation}\label{ineq8}
\Delta_{H}^{\infty}\phi(t^0,x^0)-\phi_{t}(t^0,x^0)\leq0 .
\end{equation}
Dividing (\ref{ineq6}) by $\varepsilon^2$, using (\ref{equ2})
and(\ref{ineq8}),  we get
\begin{equation}
\limsup_{\varepsilon\rightarrow
0}\frac{1}{\varepsilon^2}\left(-\phi(t^0,x^0)+\frac{1}{2}\fint_{t^0-\varepsilon^2}^{t^0}\left(\max_{y\in\overline{B}_{\varepsilon}(x^0)}\phi(s,y)
+\min_{y\in\overline{B}_{\varepsilon}(x^0)}\phi(s,y)\right)ds\right)\leq
0.
\end{equation}
That is
\begin{equation}
\phi(t^0,x^0)\geq\fint_{t^0-\varepsilon^2}^{t^0}\left(\max_{y\in\overline{B}_{\varepsilon}(x^0)}\phi(s,y)
+\min_{y\in\overline{B}_{\varepsilon}(x^0)}\phi(s,y)\right)ds+o(\varepsilon^2).
\end{equation}

Finally, let $\phi\in C_{H}^2(\Omega_{T})$ be a test function such
that $u-\phi$ has a strict minimum at $(t^0,x^0)$ and
$\nabla_{0}\phi(t^0,x^0)=0$. With the help of Theorem \ref{thm2-1},
we also assume that $(X^2\phi)^{*}(t^0,x^0)=0$, and thus the Taylor
expansion (\ref{taylor}) implies
\begin{equation*}
\phi(s,y)-\phi(t^0,x^0)=\phi_t(t^0,x^0)(s-t^0)+T\phi(t^0,x^0)\left(y_{2n+1}-x_{2n+1}^0+2\sum_{i=1}^n(y_{n+i}x_i-y_ix_{n+i})\right)
+o(\varepsilon^2).
\end{equation*}
That is
\begin{align*}
&\frac{1}{2}\left(\max_{y\in\overline{B_{\varepsilon}}(x^0)}(\phi(s,y)-\phi(t^0,x^0))+
\min_{y\in\overline{B}_{\varepsilon}(x^0)}(\phi(s,y)-\phi(t^0,x^0))\right)\\
&=\phi_t(t^0,x^0)(s-t^0)+T\phi(t^0,x^0)\{\max_{y\in\overline{B}_{\varepsilon}(x^0)}\left(y_{2n+1}-x_{2n+1}^0
+2\sum_{i=1}^n(y_{n+i}x^0_i-y_ix^0_{n+i})\right)\\
&+\min_{y\in\overline{B}_{\varepsilon}(x^0)}\left(y_{2n+1}-x^0_{2n+1}+2\sum_{i=1}^n(y_{n+i}x^0_i-y_ix^0_{n+i})\right)\}
+o(\varepsilon^2).
\end{align*}
We claim
\begin{equation}\label{claim0}
\max_{y\in\overline{B}_{\varepsilon}(x^0)}\left(y_{2n+1}-x^0_{2n+1}+2\sum_{i=1}^n(y_{n+i}x^0_i-y_ix^0_{n+i})\right)
+\min_{y\in\overline{B}_{\varepsilon}(x^0)}\left(y_{2n+1}-x^0_{2n+1}+2\sum_{i=1}^n(y_{n+i}x^0_i-y_ix^0_{n+i})\right)=0.
\end{equation}
Indeed, if $y\in \overline{B_{\varepsilon}}(x^0)$, then
\begin{align*}
&\left(y_{2n+1}-x^0_{2n+1}+2\sum_{i=1}^n(y_{n+i}x^0_i-y_ix^0_{n+i})\right)^2\\
&\leq\left(\sum\limits_{i=1}^{2n}(y_i-x^0_i)^2\right)^2+
\left(y_{2n+1}-x^0_{2n+1}+2\sum_{i=1}^n\left(y_{n+i}x^0_i-y_ix^0_{n+i}\right)\right)^2\leq
\varepsilon^4,
\end{align*}
thus
\begin{equation*}
-\varepsilon^2\leq
y_{2n+1}-x^0_{2n+1}+2\sum_{i=1}^n\left(y_{n+i}x^0_i-y_ix^0_{n+i}\right)\leq
\varepsilon^2.
\end{equation*}
Moreover, let
$$y_{max}=(x^0_1,\cdots,x^0_{2n},x^0_{2n+1}+\varepsilon^2)\in \overline{B}_{\varepsilon}(x^0),$$ and
$$y_{min}=(x^0_1,\cdots,x^0_{2n},x^0_{2n+1}-\varepsilon^2)\in \overline{B}_{\varepsilon}(x^0).$$
The maximum and minimum value can achieve at $y_{max}$ and
$y_{min}$, respectively, i.e.
\begin{equation*}
\max_{y\in\overline{B}_{\varepsilon}(x^0)}\left(y_{2n+1}-x^0_{2n+1}+2\sum_{i=1}^n(y_{n+i}x^0_i-y_ix^0_{n+i})\right)=\varepsilon^2,
\end{equation*}
and
\begin{equation*}
\min_{y\in\overline{B}_{\varepsilon}(x^0)}\left(y_{2n+1}-x^0_{2n+1}+2\sum_{i=1}^n(y_{n+i}x^0_i-y_ix^0_{n+i})\right)=-\varepsilon^2.
\end{equation*}
This ends the proof of the claim (\ref{claim0}). Therefore
\begin{equation}\label{equ3}
\frac{1}{2}\left(\max_{y\in\overline{B}_{\varepsilon}(x^0)}\left(\phi(s,y)-\phi(t^0,x^0)\right)+
\min_{y\in\overline{B}_{\varepsilon}(x^0)}\left(\phi(s,y)-\phi(t^0,x^0)\right)\right)
=\phi_t(t^0,x^0)(s-t^0)+o(\varepsilon^2).
\end{equation}

Suppose that the asymptotic mean value formula (\ref{vis2}) holds at
$(t^0,x^0)$, we get
\begin{equation*}
\phi(t^0,x^0)\geq
\frac{1}{2}\fint_{t^0-\varepsilon^2}^{t^0}\left(\max_{y\in\overline{B}_{\varepsilon}(x^0)}\phi(s,y)
+\min_{y\in\overline{B}_{\varepsilon}(x^0)}\phi(s,y)\right)ds+o(\varepsilon^2).
\end{equation*}
Hence, by (\ref{equ3}), we have
\begin{align}\label{ineq9}
0&\geq
\frac{1}{2}\fint_{t^0-\varepsilon^2}^{t^0}\left(\max_{y\in\overline{B}_{\varepsilon}(x^0)}\left(\phi(s,y)-\phi(t^0,x^0)\right)
+\min_{y\in\overline{B}_{\varepsilon}(x^0)}(\phi(s,y)-\phi(t^0,x^0))\right)ds+o(\varepsilon^2)\nonumber\\
&=\fint_{t^0-\varepsilon^2}^{t^0}\phi_t(t^0,x^0)(s-t^0)ds
+o(\varepsilon^2)\nonumber\\
&=-\frac{\varepsilon^2}{2}\phi_t(t^0,x^0)+o(\varepsilon^2).
\end{align}
Dividing (\ref{ineq9}) by $\varepsilon^2$ and passing to a limit, we
obtain
\begin{equation*}
\phi_{t}(t^0,x^0)\geq 0.
\end{equation*}
Thus, Theorem  \ref{thm2-1} shows $u$ is a viscosity super-solution
of (\ref{parabolic-infinity}).

Suppose that $u$ is a viscosity super-solution of
(\ref{parabolic-infinity}). Let $\phi\in\Omega_{T}$ be a test
function such that $u-\phi$ has a strict minimum at $(t^0,x^0)$,
$\nabla_{0}\phi(t^0,x^0)=0$ and $(X^2\phi)^{*}(t^0,x^0)=0$, we have
\begin{equation*}
\phi_{t}(t^0,x^0)\geq 0.
\end{equation*}
By (\ref{equ3}), we get
\begin{align*}
&\frac{1}{2}\fint_{t^0-\varepsilon^2}^{t^0}\left(\max_{y\in\overline{B}_{\varepsilon}(x^0)}\phi(s,y)
+\min_{y\in\overline{B}_{\varepsilon}(x^0)}\phi(s,y)\right)ds-\phi(t^0,x^0))\\
&=\fint_{t^0-\varepsilon^2}^{t^0}\phi_t(t^0,x^0)(s-t^0)ds
+o(\varepsilon^2)\\
&=-\frac{\varepsilon^2}{2}\phi_t(t^0,x^0)+o(\varepsilon^2) \leq
o(\varepsilon^2).
\end{align*}
Thus, dividing the above equality by $\varepsilon^2$ and passing to
a limit, we have
\begin{equation*}
\phi(t^0,x^0)\geq\fint_{t^0-\varepsilon^2}^{t^0}\left(\max_{y\in\overline{B}_{\varepsilon}(x^0)}\phi(s,y)
+\min_{y\in\overline{B}_{\varepsilon}(x^0)}\phi(s,y)\right)ds+o(\varepsilon^2).
\end{equation*}
Therefore, the proof of the theorem is complete. \qed
\\

Combining the case $p=2$ with the case $p=\infty$, we prove
the general case $1<p<\infty$.\\

{\textit{Proof of Theorem \ref{thm1}}}\quad Assume that $p\geq 2$ so
that $\alpha\geq 0$. Multiplying (\ref{identity}) by $\beta$,
(\ref{ineq6}) by $\alpha$ , and adding, we get
\begin{align*}
&\frac{\alpha}{2}\fint_{t-\varepsilon^2}^{t}\left(\max_{y\in\overline{B}_{\varepsilon}(x)}\phi(s,y)
+\min_{y\in\overline{B}_{\varepsilon}(x)}\phi(s,y)\right)ds
+\beta\fint_{t-\varepsilon^2}^{t}\fint_{B_{\varepsilon}(x)}\psi(x^{-1}\circ
y)\phi(s,y)dyds -\phi(x,t)\\
 &\geq
\frac{\alpha}{2}\varepsilon^2\
        \Big(\fint_{t-\varepsilon^2}^{t}\sum_{i,j=1}^{2n}X_iX_j\phi(x,t)\cdot\frac{(x_{i}^{\varepsilon,s}-x_i)}{\varepsilon}
\cdot\frac{(x_{j}^{\varepsilon,s}-x_j)}{\varepsilon}ds
-\phi_t(t,x)\Big)
+\frac{\beta}{2}\varepsilon^2\Big(M(n)\Delta_{H}\phi(t,x)-\phi_t(t,x)\Big)
+o(\varepsilon^2)\\
&=\frac{\beta}{2}\varepsilon^2\Big(\frac{\alpha}{\beta}\fint_{t-\varepsilon^2}^{t}
\sum_{i,j=1}^{2n}X_iX_j\phi(x,t)\cdot\frac{(x_{i}^{\varepsilon,s}-x_i)}{\varepsilon}
\cdot\frac{(x_{j}^{\varepsilon,s}-x_j)}{\varepsilon}ds
+M(n)\Delta_{H}\phi(t,x)
-(\frac{\alpha}{\beta}+1)\phi_t(t,x)\Big)+o(\varepsilon^2).
\end{align*}
Thanks to \begin{equation*}  \left\{
\begin{array}{lll}
\beta M(n)(p-2)=\alpha, \\[1mm]
 \alpha+\beta=1,
\end{array}
 \right.
\end{equation*}
we have
\begin{align}\label{9}
&\frac{\alpha}{2}\fint_{t-\varepsilon^2}^{t}(\max_{y\in\overline{B}_{\varepsilon}(x)}\phi(s,y)
+\min_{y\in\overline{B}_{\varepsilon}(x)}\phi(s,y))ds
+\beta\fint_{t-\varepsilon^2}^{t}\fint_{B_{\varepsilon}(x)}\psi(x^{-1}\circ
y)\phi(s,y)dyds -\phi(x,t)\nonumber\\
 &\geq
\frac{\beta}{2}\varepsilon^2\Big(M(n)(p-2)\fint_{t-\varepsilon^2}^{t}
\sum_{i,j=1}^{2n}X_iX_j\phi(x,t)\cdot\frac{(x_{i}^{\varepsilon,s}-x_i)}{\varepsilon}
\cdot\frac{(x_{j}^{\varepsilon,s}-x_j)}{\varepsilon}ds\nonumber\\
 &+M(n)\Delta_{H}\phi(t,x) -(M(n)(p-2)+1)\phi_t(t,x)\Big)+o(\varepsilon^2).
\end{align}
 The rest proof follows that of Theorem
\ref{thm-infinity}. Furthermore, by considering the maximum point
instead of the minimum point, we can get a reverse inequality to
(\ref{9}).

If $1<p<2$, it follows that $\alpha<0$, and the inequality (\ref{9})
is reversed. On the other hand, so is the reverse inequality that
can be obtained by considering the maximum point instead of the
minimum point. The argument then continues to work in the same way
as before. \qed


\end{document}